\documentclass[12pt,notitlepage,twoside]{article}

  \usepackage[all]{xy}

\usepackage{latexsym}
\usepackage{amsmath}
\usepackage{amsfonts}
\usepackage{amssymb}
\let \n = \noindent
\let \dis = \displaystyle

\newcommand{\be}{\begin{equation}}
\newcommand{\ee}{\end{equation}}
\newenvironment{pf}{\noindent{\bf Proof.}\enspace}{
\hfill$\Box$\medskip}
\newenvironment{pfn}[1]{\noindent{\bf Proof of
{#1}\enspace}}{
\hfill$\Box$\medskip}

\newtheorem{thm}{Theorem}[section]
\newtheorem{prop}[thm]{Proposition}
\newtheorem{lem}[thm]{Lemma}
\newtheorem{rem}[thm]{Remark}
\newtheorem{cor}[thm]{Corollary}
\newtheorem{defn}[thm]{Definition}
\numberwithin{equation}{section}

\author{Hichem Chtioui \footnote{ E-mail
 addresses: \texttt{Hichem.Chtioui@fss.rnu.tn}
 .}\\
{\footnotesize  Sfax University,
Faculty of Sciences of Sfax, 3018 Sfax,
Tunisia.}}

\title {On the Chen-Lin conjecture for the prescribed scalar curvature problem }

\begin{document}

\date{ }

\maketitle

{\footnotesize

\n{\bf Abstract.}  We prove a criterion of existence of solutions conjectured by C. C. Chen and C. S. Lin \cite{CL1} for the prescribed scalar curvature problem on the standard  $n$-dimensional sphere, $n\geq 3$.\\
\n{\bf  MSC 2000:}\quad    58E05, 35J60.\\
\n {\bf Key words:} Scalar curvature, Nonlinear PDEs, Variational method,  Loss  of compactness, Critical points at infinity.}




\section{Introduction}

Let $\mathbb{S}^n= \{x\in \mathbb{R}^{n+1}, |x|=1\}, n\geq 3$, endowed with the standard metric $g_0= \sum_{i=1}^{n+1} dx_i^2$ and let $K: \mathbb{S}^n\rightarrow \mathbb{R}$ be a given function. We are interested to the following scalar curvature (or Nirenberg) problem: finding suitable conditions on $K$ to be realized as the scalar curvature of some metric $g$ on $\mathbb{S}^n$ conformally equivalent to $g_0$. Writing $g= u^\frac{4}{n-2}g_0$, where $u:\mathbb{S}^n\rightarrow\mathbb{R}$ is smooth and positive, then this problem is  equivalent to solve the following nonlinear critical PDE

\begin{equation}\label{1.1}
\left\{
  \begin{array}{ccc}
    \dis -L_{g_0} u &=& \dis K   u^{\frac{n+2}{n-2}}\\
     \dis u &>& 0 \;\mbox{ on } \mathbb{S}^n.
  \end{array}
\right.
\end{equation}
Here $-L_{g_0} u= -\frac{4(n-1)}{n-2}\Delta_{g_0} u + n(n-1) u$  is the conformal Laplacian  and  $\Delta_{g_0}$ is the Laplace-Beltrami operator of $(\mathbb{S}^n, g_0)$, see \cite{A1}. The same problem can be studied on any closed -compact riemannian manifold.

 In \eqref{1.1}, the exponent on the right-hand side is $N-1$ where $N= \frac{2n}{n-2}$ is the critical case of the Sobolev embedding of $H^1(\mathbb{S}^n)$ in $L^N(\mathbb{S}^n)$. In virtue of the lack of compactness of the embedding, the variational structure associated to \eqref{1.1} does not satisfy the Palais-Smale condition, which makes the problem particulary challenging.  Moreover, besides the necessary condition that $K$ be positive somewhere there are obstructions to the existence of solutions for \eqref{1.1}, see \cite{KW1}.

\n To our knowledge, the first contribution on this problem is by D. Koutroufiotis \cite{KT1} for two-dimensions, where the analog of \eqref{1.1} has an exponential form. He has been able to solve the problem on $\mathbb{S}^2$ when $K$ is assumed to  be an antipodally symmetric function which is close to 1. Later, Moser \cite{MO1} proved that for $K$ being antipodally symmetric function which is just positive somewhere on $\mathbb{S}^2$, a solution always exists. Moser's Theorem was generalized by J. Escobar and R. Schoen \cite{ES1} for $n\geq 3$ under a suitable flatness condition up to order $n-2$. But more works dealt with nonsymmetric cases, see for example Chang-Yang \cite{CY1}and Chen-Li \cite{com}, for $n=2$, Bahri-Coron \cite{BC1} and Schoen-Zhang \cite{SZ} for $n=3$ and Chen-Lin \cite{CL1} and Li \cite{yy1} for higher dimensions, where the symmetry condition on $K$ is replaced by conditions on the critical points.

\n The scalar curvature problem continues to be one of the major topics in geometric analysis. Intensive studies have been devoted to this problem. Without attempting to give a complete list of references, one may see for example \cite{OH1}, \cite{AB}, \cite{BCH2}, \cite{BCH3}, \cite{ca}, \cite{CY2}, \cite{CGY1}, \cite{com1}, \cite{ch}, \cite{Y3}, \cite{ma2},  \cite{sho1}, \cite{sho2} and the references therein.
\vspace{0.1cm}

 In 2001, Chen-Lin \cite{CL1} considered problem \eqref{1.1} when $K$ is a $C^{1}$ positive function on $S^n$, $n\geq 3$, satisfying some flatness conditions near its critical points with a flatness order $\beta \in (1, \infty)$. Based on tricky variational tools and a refined blow-up analysis, Chen-Lin studied the blowing-up behavior of sequences of blow-up solutions of \eqref{1.1} when the prescribed scalar curvature $K$ is a small perturbation of a positive constant. More precisely, they were able to provide necessary conditions for all possible blow-up points of problem \eqref{1.1} (see Theorem 1.4 of \cite{CL1}) as well as they were able to prove an apriori bound for solutions of \eqref{1.1} under the assumption that
\begin{equation*}
    \beta(y)> \frac{n-2}{2}, \ \hbox{ for any }y \in \Gamma, \eqno{\mathbf{(A_1)}}
\end{equation*}
(see Theorems 1.1 and 10.3 of \cite{CL1}). Here $\Gamma$ denotes the set of critical points of $K$.\\
In their paper, (\cite{CL1}, Theorem 1.1), Chen-Lin conjectured a Leray-Schauder degree formula for the problem under $\mathbf{(A_1)}$ and the following condition
\begin{equation*}\label{aa}
    \dis\frac{1}{\beta^*(y_i)}+\frac{1}{\beta^*(y_j)}-\frac{2}{n-2}\neq 0, \ \hbox{for any }y_i\neq y_j \in \Gamma^{-}. \eqno{\mathbf{(A_2)}}
\end{equation*}
Here $\Gamma^-$ (as denoted in \cite{CL1}) is the set of critical points of $K$ where  blow-up phenomena may occur and $\beta^*(y)=\min(\beta(y),n)$. The predicted degree formula involves all subsets $A$ of critical points in $\Gamma^-$ such that

$$\dis\frac{1}{\beta^*(y_i)}+\frac{1}{\beta^*(y_j)}-\frac{2}{n-2}> 0, \forall y_i\neq y_j \in A.$$
Generally, there are two ways to establish the Leray-Schauder degree for problem \eqref{1.1} under the assumption that $K$ is flat near its critical points with a Morse structure. Namely, \\ \\
\n $(f)_{\beta}$: \; $K: \mathbb{S}^n\rightarrow \mathbb{R},  n\geq 3,$  is a $C^1$-positive function such that near each of its critical point $y$ there exists a real number $\beta= \beta(y)\in (1, n)$ such that in some geodesic normal coordinates system centered  at $y$, the following expansion holds

$$K(x)= K(y) +\sum_{k=1}^n b_k|(x-y)_k|^\beta+ R(x-y),$$
where   $b_k= b_k(y)\neq 0$,  $\forall k=1\ldots, n,$ $\sum_{k=1}^n b_k(y)\neq 0$ and $R$ satisfies
\begin{center}
$\sum_{j=0}^{[\beta]}|\dis\nabla^{j}\dis
R(x-y)||\dis x-y|^{j-\beta}=\dis o(1)$, as $x$ tends to $y$.
\end{center}
One way is to identify all blow-up points of problem \eqref{1.1} and compute the local degree for each blow-up solution. Another way is to identify all critical points at infinity \cite{b1} of the problem and compute the index of the associated energy functional $J$ at each critical point at infinity.\\
It is shown in (\cite{yy1}, section 6) and (\cite{L3}, Theorem 09) that under the above $(f)_{\beta}$ condition, the Leray-Schauder degree of all solutions of \eqref{1.1} is equal to an index-counting formula which involves all critical points at infinity of the problem. Note that there is a one to one correspondence of blow-up points and critical points at in infinity, see \cite{ma1} and \cite{ma3}.\\
When $\beta(y)\in (n-2,n)$ for any $y \in \Gamma$, the apriori bound of solutions and the Leray-Schauder degree formula of problem \eqref{1.1} were obtained previously by Y. Li \cite{yy1}. These results were extended to the case $\beta(y) \in [n-2,n)$ for any $y \in \Gamma$ in \cite{L3} under some additional conditions on $K$.
\vspace{0.1cm}

In the present paper we consider the case of $\beta(y) \in (1,n)$ for any $y \in \Gamma$ and our main purpose is to prove the degree counting formula conjectured by Chen-Lin \cite{CL1} and related ones in such a case.\\
We first introduce some notations and definitions extracted from \cite{CL1}. Assume that $K$ satisfies $(f)_{\beta}$-condition. Let
$$\Gamma^-=\{y\in \Gamma, \; \sum_{k=1}^n b_k(y)<0\}$$
and let $\Lambda^{-}$ be a collection of subsets of $\Gamma^{-}$ defined as follows:

\vspace{0.5cm}

\noindent {\bf Definition 1.} Let $A=\{y_1,...,y_p\}$ be a subset of $\Gamma^{-}$. $A$ is an element of $\Lambda^{-}$ if and only if $A$ satisfies the following two conditions
\begin{enumerate}
  \item $\sharp A \geq 2$
  \item $\dis\frac{1}{\beta^*(y_{i})}+\frac{1}{\beta^*(y_{j})}-\frac{2}{n-2}> 0$, for any $y_i\neq y_j$ in $A$.

\end{enumerate}

We shall prove the following Theorem which confirms the predicted Leray-Schauder degree formula announced in ( \cite{CL1}, Theorem 1.1) for $\beta \in (1,n)$.
\begin{thm}\label{th1.1}Assume $(f)_{\beta}$,  $\mathbf{(A_1)}$ and $\mathbf{(A_2)}$. Let $c$ a positive constant such that
\begin{equation*}
    \frac{1}{c}\leq w(x)\leq c, \ \forall x \in S^n,
\end{equation*}
for all solutions $w$ of \eqref{1.1}, (see \cite{CL1}, Theorem 1.1). Then for all $R\geq c$,
\begin{eqnarray*}
  d &:=& deg\big(w+L_{g_{0}}^{-1}(Kw^{\frac{n+2}{n-2}}),O_{R},0\big)\\
   &=& -\biggr[1+ \sum_{y\in \Gamma^-}(-1)^{n+1+i(y)}+ \sum_{A=\{y_1, \ldots, y_p\}\in \Lambda^-}(-1)^{p(n+1)+\sum_{j=1}^p  i(y_j)}\biggr],
\end{eqnarray*}
where $deg$ denotes the Leray-Schauder degree in $C^{2,\alpha}(S^n)$, $0<\alpha<1$, $O_R$ is defined by
\begin{equation*}
O_R =\big\{v \in C^{2,\alpha}(S^n), \ 0<\alpha<1, \ s.t, \ \frac{1}{R}\leq v \leq R \hbox{ and }\|v\|_{ C^{2,\alpha}}\leq R\big\}
\end{equation*}
and
$$i(y)= \sharp \{b_k(y), 1\leq k\leq n, s.t., b_k(y)<0\}.$$
\end{thm}
\textbf{Remarks}\\
\begin{enumerate}
  \item The above degree formula $d$ is identical to the one of Chen-Lin, since under $(f)_{\beta}$-condition the local degree $deg_{Loc}(\nabla K,y_j)=(-1)^{i(y_j)}$.
  \item If $\beta(y)\in (n-2,\infty)$ for any $y \in \Gamma$, then
  $$\dis\frac{1}{\beta^*(y_{i})}+\frac{1}{\beta^*(y_{j})}-\frac{2}{n-2}< 0, \ \forall y_i\neq y_j \in \Gamma. $$
  In this case $\Lambda^{-}$ is empty. In particular if $\beta(y)\in (n-2,n)$ for any $y \in \Gamma$, the degree counting formula of Theorem \ref{th1.2} is reduced to
  $$d=-1+ \sum_{y\in \Gamma^-}(-1)^{n+i(y)}.$$
  It is the Leray-Schauder degree formula proved in (\cite{yy1}, Theorem 0.1).
  \item The above Theorem can be extended to the case of $\beta(y)\in(1,n]$ with a minor additional arguments in our proof. But its extension to the case of $\beta(y)\in(1,\infty)$ still remains open.
\end{enumerate}
\vspace{0.1cm}

In the next, we are mainly concerned with what happens to the Leray-Schauder degree formula of problem \eqref{1.1} if $K$ does not satisfy condition $\mathbf{(A_2)}$. In (\cite{CL1}, Theorem 10.3) it is shown that the apriori bound of solutions remains when $\mathbf{(A_2)}$ is removed and replaced by conditions $\mathbf{(H_1)}$ and $\mathbf{(H_2)}$ below. But no information about the degree counting formula was delivered.\\
The aim of the next Theorem is to extend the degree formula of Theorem \ref{th1.1} to functions $K$ failing $\mathbf{(A_2)}$ condition and provide an entirely new one. We need to state more notations and assumptions. \\
Let $G$ be the Green function for the operator $-L_{g_0}$ on $\mathbb{S}^n$. For each p-tuple, $p\geq 2$, of distinct critical points $\tau_p=(y_1, \ldots, y_p)\in (\Gamma^-)^p$  such that $\beta(y_i)= n-2, \forall i=1, \ldots, p$, we define a symmetric matrix $M(\tau_p)= (m_{ij})_{1\leq i\neq j \leq p}$ by:

$$m_{ij}= m(y_i, y_j)= - \tilde{c}_1 2^\frac{n-2}{2} \frac{G(y_i, y_j)}{\Big(K(y_i)K(y_j)\Big)^\frac{n-2}{4}}, \;\;\forall 1\leq i\neq j \leq p,$$and
$$m_{ii}= m(y_i, y_i)= - c_i  \frac{\sum_{k=1}^n b_k(y_i)}{K(y_i)^\frac{n}{2}}, \;\;\forall 1\leq i \leq p,$$
where $\tilde{c}_1= \dis\int_{\mathbb{R}^n}\frac{dx}{(1+|x|^2)^\frac{n+2}{2}}$ and $c_i= \dis \int_{\mathbb{R}^n}\frac{|x_1|^{\beta(y_i)}dx}{(1+|x|^2)^n}$.

\n In $c_i$, $x_1$ denotes the first component of $x$ in some geodesic normal coordinates system.\\

\n $\mathbf{(H_1)}$  Assume that for any $\tau_p =(y_1, \ldots, y_p)\in (\Gamma^-)^p$, such that $\beta(y_i)= n-2, \forall i=1, \ldots, p$ and  $y_i\neq y_j$, $\forall 1\leq i\neq j\leq p$, the least eigenvalue $\rho(\tau_p)$ of $M(\tau_p)$ is non zero.

\vspace{0.1cm}

\n For any critical point $y_{j_0}$ of $K$ such that $\beta(y_{j_0})>n-2$, we define $$B_{j_0}= \{y\in \Gamma, s.t., \dis\frac{1}{\beta(y)}+\frac{1}{\beta(y_{j_0})}-\frac{2}{n-2}= 0\}.$$

\n $\mathbf{(H_2)}$  For any p-tuple, $p\geq 1$, of distinct critical points $\tau_p =(z_1, \ldots, z_p)\in B_{j_0}^p$, we assume that

$$\tilde{\rho}(\tau_p)= \sum_{i=1}^p c_i \frac{|\sum_{k=1}^n b_k(z_i)|}{K(z_i)^\frac{n}{2}}\bigg(\frac{\tilde{c}_1 2^\frac{n-2}{2}}{c_i} \frac{K(z_i)^\frac{n}{2}}{|\sum_{k=1}^n b_k(z_i)|}\frac{G(z_i, y_{j_0})}{(K(z_i) K(y_{j_0}))^\frac{n-2}{4}}\bigg)^\frac{2\beta(y_{j_0})}{n-2}$$$$+ c_{j_0} \frac{\sum_{k=1}^n b_k(y_{j_0})}{K(y_{j_0})^\frac{n}{2}}\neq 0.$$
We now introduce the notation of $\tilde{\Lambda}^-$. Let $\tilde{\Lambda}^-$ be a collection of subsets of $\Gamma^-$ defined as follows:

\vspace{0.5cm}
\noindent {\bf Definition 2.}
Let $A=\{y_1,...,y_p\}$ be a subset of $\Gamma^-$. $A$ is an element of $\tilde{\Lambda}^-$ if and only if $A$ satisfies the following three conditions
\begin{enumerate}
  \item $\sharp A \geq 2$
  \item  $\dis\frac{1}{\beta^*(y_{i})}+\frac{1}{\beta^*(y_{j})}-\frac{2}{n-2}\geq 0$, for any $y_i \neq y_j \in A$.
  \item If $y_{\ell_1}, \ldots, y_{\ell_q}$ are all elements of $A$ such that for any $i=1,\ldots,q$ there exists $j\neq i$, $j=1,\ldots,q$ satisfying $$\dis\frac{1}{\beta^*(y_{i})}+\frac{1}{\beta^*(y_{j})}-\frac{2}{n-2}=0,$$ then $\beta(y_{\ell_1})= \ldots = \beta(y_{\ell_q})=n-2$ and $\rho(y_{\ell_1}, \ldots, y_{\ell_q})>0$.
\end{enumerate}

\begin{thm}\label{th1.2}
Assume $(f)_{\beta}$, $\mathbf{(A_1)}$, $\mathbf{(H_1)}$ and $\mathbf{(H_2)}$. Let $c$ a positive constant such that
\begin{equation*}
    \frac{1}{c}\leq w(x)\leq c, \ \forall x \in S^n,
\end{equation*}
for all solutions $w$ of \eqref{1.1}, (see \cite{CL1}, Theorem 10.3). Then for all $R\geq c$,
\begin{eqnarray*}
  d &:=& deg\big(w+L_{g_{0}}^{-1}(Kw^{\frac{n+2}{n-2}}),O_{R},0\big)\\
   &=& -\biggr[1+ \sum_{y\in \Gamma^-}(-1)^{n+1+i(y)}+ \sum_{A=\{y_1, \ldots, y_p\}\in \tilde{\Lambda}^-}(-1)^{p(n+1)+\sum_{j=1}^p  i(y_j)}\biggr].
\end{eqnarray*}
\end{thm}
\vspace{0.1cm}

Our used method to prove Theorems \ref{th1.1} and \ref{th1.2}  hinges on the critical points at infinity theory of A. Bahri \cite{b1}. We completely describe the critical points at infinity of the problem under conditions $\mathbf{(A_1)}$, $\mathbf{(H_1)}$ and $\mathbf{(H_2)}$. While condition $\mathbf{(A_1)}$ was used in \cite{CL1} as a necessary condition to prove an apriori bound of solutions of \eqref{1.1}, $\mathbf{(A_1)}$ is used here as a necessary condition to prove that there is no critical point at infinity constructed by a solution of \eqref{1.1} and a sum of bubbles (solutions of \eqref{1.1} when K=1).\\
\vspace{0.1cm}
The computation of the Leray-Schauder degree for problem \eqref{1.1} is useful. For instance we illustrate its usefulness through the following existence and multiplicity result.
\begin{thm}\label{th1.3}
Under assumptions $(f)_{\beta}$, $\mathbf{(A_1)}$, $\mathbf{(H_1)}$ and $\mathbf{(H_2)}$, equation \eqref{1.1} admits a solution provided $d\neq 0$. Moreover, for generic $K$, the number of solutions of \eqref{1.1} is larger or equal to $\mid d \mid$.
\end{thm}
\vspace{0.1cm}

The last result of this paper is devoted to establish a general index-counting formula of existence of solutions extending the ones of the above Theorems to the case of functions $K$ failing $\mathbf{(A_1)}$ assumption. We shall prove the following
\begin{thm}\label{th1.4}
Under assumptions $(f)_{\beta}$, $\mathbf{(H_1)}$ and $\mathbf{(H_2)}$, if
$$1+\sum_{y\in \Gamma^-}(-1)^{n+1+i(y)}+ \sum_{A=\{y_1, \ldots, y_p\}\in \tilde{\Lambda}^-}(-1)^{p(n+1)+\sum_{j=1}^p  i(y_j)}\neq 0,$$
then \eqref{1.1} admits a solution.
\end{thm}
\textbf{Remark}\\
We emphasis that a necessary condition for the above index formula to be non zero is that any local maximum $y_0$ of $K$ have $\beta(y_0)>\frac{n-2}{2}$.

\vspace{0.1cm}

The proof of our results will be performed in Section 4. In Section 3,  Theorem \ref{th3.3} completely describes when and where critical points at infinity occur under the assumption that \eqref{1.1} has no solution. In particular we will show that the possibility of tower bubbles phenomenon is excluded. The result is valid for any flatness order $\beta$ in $(1,n)$. It is interesting in its self. Note that assumption $\mathbf{(A_1)}$ is optimal in the used method of \cite{CL1}, since the apriori bound Theorem fails for $\beta<\frac{n-2}{2}$, see \cite{CL11}. \\
In the next section we set up the associated variational structure and we recall some preliminary results.\\

\n \textbf{Acknowledgement:} Part of this work was done when the author enjoys the hospitality of Giessen University. He takes the opportunity to acknowledge the excellent working condition and to thank Professor Mohameden Ould Ahmedou for many discussions about the subject of this paper.

\section{Preliminaries}
The Sobolev space $H^1(\mathbb{S}^n)$ is equipped with the norm $\dis \|u\|= \Big(\int_{\mathbb{S}^n} -L_{g_0} u u dv_{g_0} \Big)^\frac{1}{2}$. Let
$$\Sigma=\{u\in H^1(\mathbb{S}^n), \|u\|=1\} \mbox{ and } \Sigma^+=\{u\in \Sigma, u>0\}.$$
Up to a multiplicative constant, a solution of \eqref{1.1} is a critical point of the variational functional
$$J(u)= \frac{\|u\|^2}{\bigg(\int_{\mathbb{S}^n}K u^\frac{2n}{n-2}\bigg)^\frac{n-2}{2}}, \; u\in \Sigma^+.$$
Since the Sobolev-embedding $H^1(\mathbb{S}^n)\hookrightarrow L^\frac{2n}{n-2}(\mathbb{S}^n)$ is not compact, the functional $J$ does not satisfy the Palis-Smale condition ((P.S) for short). Its failure can be descried as follows:

\n \n For any $a\in \mathbb{S}^n$ and $\lambda>0$, we denote

$$\delta_{(a, \lambda)}(x)= c_0\bigg(\frac{\lambda}{\lambda^2+1+(1-\lambda^2)\cos d(a, x)}\bigg)^\frac{n-2}{2}, x\in \mathbb{S}^n.$$
Here $c_0$ is a fixed constant chosen so that $\delta_{(a, \lambda)}$ is a solution of the Yamabe problem
$$-L_{g_0} u = u^\frac{n+2}{n-2}, \quad  u>0 \mbox{ on } \mathbb{S}^n.$$

\n For $p\in \mathbb{N}^*$ and $\varepsilon>0$, we set
\begin{eqnarray*}
V(p, \varepsilon)&=& \biggr\{u\in \Sigma^+, \exists a_1, \ldots, a_p\in \mathbb{S}^n,  \exists \lambda_1, \ldots, \lambda_p > \varepsilon^{-1}, \exists \alpha_1, \ldots, \alpha_p>0 \mbox{ satsifying }\\
 &&\|u-\sum_{i=1}^p \alpha_i \delta_{(a_i, \lambda_i)}\|<\varepsilon, \ \big|\frac{ \alpha_i^\frac{4}{n-2} K(a_i)} {\alpha_j^\frac{4}{n-2} K(a_j)} -1\big|<\varepsilon,\mbox{ and } \varepsilon_{ij}< \varepsilon,  \forall 1\leq i\neq j \leq p\biggr\}.
 \end{eqnarray*}
 Here $\dis \varepsilon_{ij}= \bigg(\frac{\lambda_i}{\lambda_j} + \frac{\lambda_j}{\lambda_i} + \frac{\lambda_i \lambda_j}{2} (1-\cos d(a_i, a_j))\bigg)^\frac{2-n}{2}$. For  $w$ a solution of \eqref{1.1} or zero, we also define\\
\begin{eqnarray*}
V(p,\varepsilon,w)&=& \biggr\{u\in \Sigma^+, \exists a_1, \ldots, a_p\in \mathbb{S}^n,  \exists \lambda_1, \ldots, \lambda_p > \varepsilon^{-1}, \exists \alpha_0,\alpha_1, \ldots, \alpha_p>0\\
&& \mbox{ satsifying }\|u-(\alpha_0 w+\sum_{i=1}^p \alpha_i \delta_{(a_i, \lambda_i)})\|<\varepsilon, \ \big|\frac{ \alpha_i^\frac{4}{n-2} K(a_i)} {\alpha_j^\frac{4}{n-2} K(a_j)} -1\big|<\varepsilon,\\
&&\varepsilon_{ij}< \varepsilon,  \forall 1\leq i\neq j \leq p\mbox{ and }\big|\displaystyle\frac{\alpha_{0}^{\frac{4}{n-2}}}{\alpha_i^\frac{4}{n-2} K(a_i)}-1\big|<\varepsilon, \ \forall 1\leq i \leq p\biggr\}.
\end{eqnarray*}
Of course $V(p,\varepsilon,w)=V(p,\varepsilon)$ if $w=0$.
\begin{prop}\label{prop1} \cite{b2}, \cite{bc2}, \cite{S3} Let $(u_k)$ be a sequence in $\Sigma^+$ such that
$J(u_k)$ is bounded and $\partial J(u_k)$ goes to zero. Then there
exist an integer $p \in \mathbb{N}^*$, a sequence $(\varepsilon_k)
>0$, $\varepsilon_k$ tends to zero, and an extracted subsequence
of $u_k$'s, again denoted $(u_k)$, such that $u_k \in
V(p,\varepsilon_k ,w)$ where $w$ is zero or a solution of
\eqref{1.1}.
\end{prop}

\n Following \cite{b2}, we have the following parametrization of $V(p, \varepsilon,w)$. Let $W_{s}(w)$ and $W_{u}(w)$ denote (respectively) the stable and unstable manifolds for a decreasing pseudo-gradient of $J$ at $w$.
\begin{prop}\label{prop2}
 For any $p \in \mathbb{N}^*$, there is $\varepsilon_{p}>0$ such that if
$\varepsilon\leq\varepsilon_{p}$ and $u\in V(p,\varepsilon,w)$, then
the following minimization
problem$$\min_{\begin{array}{cc}\alpha_{i}>0,\lambda_{i}>0,a_{i}\in
S^{n},\\h\in T_{w}(W_{u}(w)). \end{array}}\Big\| u
-\dis\sum_{i=1}^{p}
\alpha_{i}\delta_{(a_{i},\lambda_{i})}-\alpha_{0}(w+h)\Big\|,$$ has
a unique solution $(\alpha,\lambda,a,h)$, up to a permutation.
\end{prop}
In particular, we can write $u$ as follows
$$u=\dis\sum_{i=1}^{p}\alpha_{i}\delta_{(a_{i},
\lambda_{i})}+\alpha_{0} (w+h)+v,$$ where $v$ belongs to
$H^{1}(S^{n})\cap T_{w}(W_{s}(w))$ and satisfies the following:\\
$$
(V_{0}): \begin{cases}
    \Big\langle v ,\psi \Big\rangle=0 \; for\;
    \psi\in\{\displaystyle\delta_{i},\frac{\displaystyle\partial\delta_{i}}{\displaystyle\partial\lambda_{i}},\frac{\displaystyle\partial\delta_{i}}
    {\displaystyle\partial
    a_{i}},i=1,...,p\}\\\\
    \Big\langle v,w\Big\rangle=0\\\\
    \Big\langle v,h\Big\rangle=0 \; for \;all\; h\in T_{w}W_{u}(w)
  \end{cases}$$ here, $\delta_i
=\delta_{(a_{i},\lambda_{i})}$ and $<.,.>$ denotes the scalar
product defined on $H^1(S^n)$ by
$$<w, w'>= \int_{\mathbb{S}^n}-L_{g_0} w w' dv_{g_0},  \; \forall w, w' \in H^1(\mathbb{S}^n).$$
Next, we say that $v\in (V_0)$ if $v$ satisfies $(V_0)$.

\begin{prop}\label{} \cite{b2}, \cite{bc2} There is a $\mathcal{C}^{1}$-map which to each
$(\alpha_{i},a_{i},\lambda_{i},h)$ such that\\ $\dis\sum_{i=1}^{p}
\alpha_{i}\delta_{(a_{i},\lambda_{i})}+\alpha_{0}(w+h)$ belongs to
$V(p,\varepsilon,w)$ associates $\overline{v}=\overline
v(\alpha,a,\lambda,h)$ such that $\overline{v}$ is unique and
satisfies:
$$J\Big(\sum_{i=1}^{p}
\alpha_{i}\delta_{(a_{i},\lambda_{i})}+\alpha_{0}(w+h)+\overline{v}\Big)=\dis\min_{v
\in (V_0)}\Big\{J\Big(\sum_{i=1}^{p}
\alpha_{i}\delta_{(a_{i},\lambda_{i})}+\alpha_{0}(w+h)+v\Big)\Big\}.$$
Moreover,there exists a change of variables
$v-\overline{v}\rightarrow V$ such that $$J\Big(\dis\sum_{i=1}^{p}
\alpha_{i}\delta_{(a_{i},\lambda_{i})}+\alpha_{0}(w+h)+v\Big)=J\Big(\dis\sum_{i=1}^{p}
\alpha_{i}\delta_{(a_{i},\lambda_{i})}+\alpha_{0}(w+h)+\overline{v}\Big)+\parallel
V \parallel^{2}.$$
\end{prop}

\n The estimate of $\|\bar{v}\|$ was provided in \cite{cbd1} under $(f)_\beta$-condition, $1< \beta< n$.

\begin{prop}\label{prop3}\cite{cbd1} Under $(f)_\beta$-condition, $\beta\in (1, n)$, there exists $c>0$
such that
$$\|\bar{v}\|\leq c \sum_{i=1}^p \bigg(\frac{1}{\lambda_i^\frac{n}{2}} + \frac{1}{\lambda_i^\beta}+ \frac{|\nabla K(a_i)|}{\lambda_i}+ \frac{(\log \lambda_i)^\frac{n+2}{2n}}{\lambda_i^\frac{n+2}{2}}\bigg)+c\left\{
                                                                 \begin{array}{ll}
                                                                   \sum_{k\neq r}\varepsilon_{kr}^\frac{n+2}{2(n-2)}(\log \varepsilon_{kr}^{-1})^\frac{n+2}{2n}, & \hbox{ if } n\geq 6, \\
                                                                   \sum_{k\neq r}\varepsilon_{kr}(\log \varepsilon_{kr}^{-1})^\frac{n-2}{n}, & \hbox{ if } n\leq 5.
                                                                 \end{array}
                                                               \right.
$$
\end{prop}

\n We now state the definition of a critical point at infinity.

\begin{defn}\cite{b1}, \cite{b2} A critical point at infinity of $J$on $\Sigma^{+}$ is a limit of a
flow line $u(s)$ of the equation
$$\begin{cases}
\frac{\dis\partial u}{\dis\partial s}=-\partial J(u(s))\\
u(0)=u_{0}
\end{cases}$$
such that $u(s)$ remains in $V(p,\varepsilon(s),w)$ for $s\geq
s_{0}$.  Here $w$ is either zero or a solution of \eqref{1.1} and
$\varepsilon(s)$ is some positive  function tending to zero when
$s\rightarrow+\infty$. Using proposition \ref{prop2}, $u(s)$ can be
written as:
\begin{center}
$u(s)=\dis \sum_{i=1}^{p}\dis\alpha_{i}(s)\dis\delta_{\dis(\dis
a_{i}(s),\dis\lambda_{i}(s))}+\dis\alpha_{0}(s)(w+h(s))+v(s)$.
\end{center}
Denoting $\dis\widetilde{\alpha}_{i}:=\dis \lim_{s \longrightarrow
+\infty} \dis\alpha_{i}(s)\;\; \mbox{,}\; \widetilde{y}_{i}:=\lim_{s
\longrightarrow + \infty} \dis a_{i}(s)$, we denote by
\begin{center}
$ \dis
\sum_{i=1}^{p}\dis\widetilde\alpha_{i}\dis\delta_{\dis(\dis\widetilde
y_{i},\infty)}+\dis\widetilde\alpha_{0}w$ or $(\dis\widetilde
y_{1},...,\dis\widetilde y_{p},w)_{\infty}$
\end{center}
 such a critical point at
infinity.
\end{defn}

\n We conclude this section by providing the expansions of $\dis \lambda_i \frac{\partial J(u)}{\partial \lambda_i}$ and $\dis \frac{1}{\lambda_i}\frac{\partial J(u)}{\partial a_i}$ in $V(p, \varepsilon)$, $p\in \mathbb{N}^*$. These expansions can be found in ( \cite{kh5}, Propositions 3.5 and 3.6 ). Let $\rho_0$ be a fixed  positive constant small enough.

\begin{prop} \label{prop4}  Let $u= \dis \sum_{i=1}^p \alpha_i \delta_{(a_i, \lambda_i)}\in V(p, \varepsilon)$ such that $a_i\in B(y_i, \rho_0), y_i\in \Gamma, \forall i=1, \ldots, p$. Under $(f)_\beta$-condition, $\beta\in (1, n)$, we have the following four expansion:

$$(i)\quad <\partial J(u), \alpha_i \lambda_i \frac{\partial \delta_{(a_i,\lambda_i)}}{\partial \lambda_i}>= -2J(u) \bigg[\tilde{c_1} \sum_{j\neq i} \alpha_i \alpha_j \lambda_i \frac{\partial \varepsilon_{ij}}{\partial \lambda_i} - \frac{n-2}{2} \alpha_i^2 \frac{c_i}{K(a_i)}\frac{\sum_{k=1}^n b_k(y_i)}{\lambda_i^{\beta(y_i)}}\bigg]$$$$+O\bigg(|a_i-y_i|^{\beta(y_i)}\bigg)+O\bigg(\frac{|a_i-y_i|^{\beta(y_i)-1}}{\lambda_i}\bigg)+ o\Big(\frac{1}{\lambda_i^{\beta(y_i)}}\Big)+ o\Big(\sum_{j\neq i}\varepsilon_{ij}\Big).$$

$$(ii)\quad <\partial J(u), \alpha_i \lambda_i \frac{\partial \delta_{(a_i,\lambda_i)}}{\partial \lambda_i}>= -2J(u) \tilde{c_1} \sum_{j\neq i} \alpha_i \alpha_j \lambda_i \frac{\partial \varepsilon_{ij}}{\partial \lambda_i} + \Big(\sum_{j=2}^{[\beta(y_i)]} \frac{|a_i-y_i|^{\beta(y_i)-j}}{\lambda_i^{j}}\Big)$$$$+O\Big(\frac{1}{\lambda_i^{\beta(y_i)}}\Big)+ o\Big(\sum_{j\neq i}\varepsilon_{ij}\Big).$$
For $k=1, \ldots, n$, we denote $(a_i)_k$ the $k^{th}$ component of $a_i$ in some geodesic normal coordinates system. We then have
$$(i)'\quad <\partial J(u), \alpha_i \frac{1}{\lambda_i} \frac{\partial \delta_{(a_i,\lambda_i)}}{\partial (a_i)_k}>= - \tilde{c}_2 \alpha_i^2 \frac{J(u)}{K(a_i)}\beta(y_i) \mbox{ sign } (a_i-y_i)_k |(a_i-y_i)_k|^{\beta(y_i)-1}$$$$ + O \bigg(\sum_{j=2}^{[\beta(y_i)]} \frac{|a_i-y_i|^{\beta(y_i)-j}}{\lambda_i^j}\bigg)+ O\Big(\frac{1}{\lambda_i^{\beta(y_i)}}\Big)+ O\Big(\sum_{j\neq i}\frac{1}{\lambda_i}\Big|\frac{\partial \varepsilon_{ij}}{\partial a_i}\Big|\Big),$$
where $\tilde{c}_2= \int_{\mathbb{R}^n}\frac{dx}{(1+|x|^2)^n}$. Moreover, if $\lambda_i|a_i-y_i|$ is bounded, then
$$(ii)'\quad <\partial J(u), \alpha_i \frac{1}{\lambda_i} \frac{\partial \delta_{(a_i,\lambda_i)}}{\partial (a_i)_k}>= -(n-2) \alpha_i^2 \frac{J(u)}{K(a_i)}\frac{b_k(y_i)}{\lambda_i^{\beta(y_i)}}\int_{\mathbb{R}^n} \frac{x_k |x_k+ \lambda_i(a_i-y_i)_k|^{\beta(y_i)}}{(1+|x|^2)^{n+1}}dx$$$$+ O\Big(\sum_{j\neq i}\frac{1}{\lambda_i}\Big|\frac{\partial \varepsilon_{ij}}{\partial a_i}\Big|\Big).$$
\end{prop}

\section{Critical points at infinity}
In this section we characterize the critical points at infinity of  problem \eqref{1.1}. Such a characterization is of course quite difficult and technical. In principle, it relies on the construction of a special pseudogradient $\widetilde{W}$ at infinity, as in (\cite{AB}, \cite{bc2}, \cite{cbd1}). The construction is more difficult when the function $K$ satisfies $(f)_{\beta}$-condition with flatness order $\beta(y)$ varies in $(1, n)$ for any $y\in \Gamma$ and leads to a new interesting phenomenon drastically different from the previous ones. Indeed, when studying the possible formations of blow-up points for the problem, it turns out that the mutual interaction among two different bubbles dominates the self interaction of the bubbles if $\beta(y)\in (n-2, n)$ for any $y\in \Gamma$. This causes all blow-up points to be single \cite{yy1}. If $\beta(y)\in (1, n-2)$ for any $y\in \Gamma$, the reverse phenomenon happens \cite{cbd1}. While if $\beta(y)= n-2$ for any $y\in \Gamma$, we have a balance phenomenon, \cite{cbd1}.\\
Now, if $\beta$ varies in $(1, n)$, particulary if we single out two bubbles at two critical points $y_i$ and $y_j$ of $K$ with $\beta(y_i)< n-2$ and $\beta(y_j)\geq n-2$, the above three phenomena may occur and each phenomenon will be related to the sign of
$$\dis\frac{1}{\beta(y_i)}+\dis\frac{1}{\beta(y_j)}-\frac{2}{n-2}.$$
Note that the
sign of $\dis\frac{1}{\beta(y_i)}+\frac{1}{\beta(y_j)}-\frac{2}{n-2}$ cannot be changed only if $\beta(y_i)$ varies in $(1, n-2)$ and  $\beta(y_j)$ varies in $(n-2, n)$.  Therefore, the aspect of the problem in the present situation is not classic and the question related to identify the the critical points at infinity of the problem will be difficult and extremely technical.
\begin{defn}\cite{b1} Let $W$ be a given vector field in $V(p,\varepsilon), \ p\geq 1$. We say that the Palais-Smale condition is satisfied along the flow lines of  $W$, or simply $W$ satisfies (P.S) condition, if $\max_{1\leq i \leq p}(\lambda_i(s))$ is a bounded function for any flow line $u(s)=\dis\sum_{i=1}^{p}
\alpha_{i}(s)\delta_{(a_{i}(s),\lambda_{i}(s))}$ of  $W$.
\end{defn}
The following Theorem shows that for any $p\geq 1$, there exists a suitable pseudogradient of $J$ in $V(p, \varepsilon)$, for which (P.S) condition is satisfied along its flow lines as long as these flow lines do not enter in neighborhood of critical points $y_i, i=1, \ldots, p$ such that $y_1\in \Gamma^{-}$ if $p=1$ and $y_i\neq y_j, \forall 1\leq  i\neq j\leq p$ with $\{y_1, \ldots, y_p\}\in  \tilde{\Lambda}^-$, if $p\geq 2$.
\begin{thm}\label{th3.1}
Let $\beta=\max\{\beta(y), y\in \Gamma\}$. Under assumptions  $(f)_{\beta}$, $\mathbf{(H_1)}$ and $\mathbf{(H_2)}$, there exist a bounded pseudogradient $\widetilde{W}$ in $V(p, \varepsilon), p\geq 1$, and a fixed positive constant $c$ such that for any $u= \dis \sum_{i=1}^p \alpha_i \delta_{(a_i, \lambda_i)}\in V(p, \varepsilon)$, the following hold
$$
(a) \; <\partial J(u), \widetilde{W}(u)> \leq -c \bigg(\sum_{i=1}^p\Big(\frac{1}{\lambda_i^{\beta}}+ \frac{|\nabla K(a_i)|}{\lambda_i}\Big)+ \sum_{j\neq i}\varepsilon_{ij}\bigg),$$
$$
(b)\;  <\partial J(u + \bar{v}), \widetilde{W}(u) + \frac{\partial \bar{v}}{\partial(\alpha_i, a_i, \lambda_i)}(\widetilde{W}(u))> \leq -c \bigg(\sum_{i=1}^p\Big(\frac{1}{\lambda_i^{\beta}}+ \frac{|\nabla K(a_i)|}{\lambda_i}\Big)+ \sum_{j\neq i}\varepsilon_{ij}\bigg).$$
\n$(c)$ \; All the component  $\lambda_i(s), \; i=1, \ldots, p,$ remain bounded as long as the associated flow line $u(s)= \sum_{i=1}^p \alpha_i(s) \delta_{(a_i(s), \lambda_i(s))}$ is out side a small neighborhood $\mathcal{N}(p, \varepsilon)$ of $\dis \sum_{i=1}^p \frac{1}{K(y_i)^\frac{n-2}{2}}\delta_{(y_i, \infty)}$, such that $y_1\in \Gamma^{-}$ if $p=1$ and $y_i\neq y_j, \forall 1\leq  i\neq j\leq p$ with $\{y_1, \ldots, y_p\}\in  \tilde{\Lambda}^-$, if $p\geq 2$. However, if $u(s)$ enter $\mathcal{N}(p, \varepsilon)$, all concentration $\lambda_i(s), i=1, \ldots, p$ tend to $\infty$.
\end{thm}
\vspace{0.1cm}

For the whole next construction, we make use the following notations and definitions.\\
Let $u=\dis\sum_{i=1}^{p}
\alpha_{i}\delta_{(a_{i},\lambda_{i})} \ \in \ V(p,\varepsilon), \ p\geq 1,$ such that $a_i \in B(y_i,\rho_{0})$, $y_i \in \Gamma$, $\forall i=1,...,p$.\\
For any index $i=1,...,p$, we define the following vector fields.
\begin{equation}\label{2.A}
    Z_i(u)=\alpha_{i}\lambda_{i}\dis \frac{ \partial \delta_{(a_{i},\lambda_{i})}}{\partial \lambda_{i}},
\end{equation}
and
\begin{equation*}
    Y_i(u)=\varphi_{M}\big(\lambda_i|a_i-y_i|\big)\frac{\alpha_i}{\lambda_i}\dis \frac{ \partial \delta_{(a_{i},\lambda_{i})}}{\partial a_{i}}\dis \frac{\nabla K(a_{i})}{|\nabla K(a_{i})|},
\end{equation*}
where $\varphi_{M}, \ M>0$, is the cut-off function defined by $\varphi_{M}(t)=0$, if $|t|\leq M$ and $\varphi_{M}(t)=1$, if $|t|\geq 2M$.\\
Observe that under the action of $Z_i$, $\lambda_i$ moves according to the differential equation
\begin{equation*}
    \dot{\lambda}_i=\lambda_i
\end{equation*}
and under the action of $Y_i$, $a_i$ moves according to
\begin{equation*}
    \dot{a}_i=\dis \frac{\nabla K(a_{i})}{|\nabla K(a_{i})|}, \ \hbox{if }a_i\neq y_i.
\end{equation*}
Moreover under $(f)_\beta$-condition we have
\begin{equation*}
    \frac{\partial K}{\partial x_k}(a_i)\sim \beta_i b_k \mbox{ sign } (a_i-y_i)_k |(a_i-y_i)_k|^{\beta_i-1}
\end{equation*}
and
\begin{equation*}
    \dis |\nabla K(a_i)|\sim \beta_i|a_i-y_i|^{\beta_i-1}.
\end{equation*}
Therefore, up to a positive multiplicative constant
\begin{equation}\label{2.B}
    Y_i(u)=\varphi_{M}\big(\lambda_i|a_i-y_i|\big)\frac{\alpha_i}{\lambda_i}\dis\sum_{k=1}^{n}b_k\mbox{ sign } (a_i-y_i)_k \frac{|(a_i-y_i)_k|^{\beta_i-1}}{|a_i-y_i|^{\beta_i-1}}\frac{ \partial \delta_{(a_{i},\lambda_{i})}}{\partial (a_{i})_k}.
\end{equation}

\begin{pfn}{\textbf{Theorem \ref{th3.1}}}
Let $\delta>0$ small enough.  We decompose $V(p, \varepsilon)$ on four  parts.

$$V_1(p, \varepsilon)= \{u=\dis\sum_{i=1}^p \alpha_i \delta_{(a_i, \lambda_i)}\in V(p, \varepsilon), s.t., a_i\in B(y_i, \rho_0), y_i\in \Gamma, \lambda_i|a_i-y_i|<\delta,$$$$ \forall i=1, \ldots, p \mbox{ and } y_i\neq y_j, \forall 1\leq i\neq j \leq p\}.$$
$$V_2(p, \varepsilon)= \{u=\dis\sum_{i=1}^p \alpha_i \delta_{(a_i, \lambda_i)}\in V(p, \varepsilon), s.t., a_i\in B(y_i, \rho_0), y_i\in \Gamma,  \forall i=1, \ldots, p,  \; y_i\neq y_j,$$$$ \forall 1\leq i\neq j \leq p, \mbox{ and there exists at least $i_0$ such that }  \lambda_{i_0} |a_{i_0}-y_{i_0}|>\frac{\delta}{2}\}.$$
$$V_3(p, \varepsilon)= \{u=\dis\sum_{i=1}^p \alpha_i \delta_{(a_i, \lambda_i)}\in V(p, \varepsilon), s.t., a_i\in B(y_i, \rho_0), y_i\in \Gamma, \forall i=1, \ldots, p$$$$ \mbox{ and there exist } i_0\neq j_0 \mbox{ such that } y_{i_0}=y_{j_0} \}.$$
$$V_4(p, \varepsilon)= \{u=\dis\sum_{i=1}^p \alpha_i \delta_{(a_i, \lambda_i)}\in V(p, \varepsilon), s.t.,\mbox{ there exists } i, \ 1\leq i \leq p,$$
$$ \mbox{ satisfying } a_i \not\in B(y,\rho_0), \ \forall y\in\Gamma \}.$$
We construct in each region $V_i(p, \varepsilon), i=1, 2, 3,4$ an appropriate pseudogradient $\widetilde{W}_i$ and then glue up through a convex combination.\\

\n $\bullet$ \underline{Pseudogradient in $V_1(p, \varepsilon), p=1$:}
\vskip0.15cm
\n Let $u= \alpha_1 \delta_{(a_1, \lambda_1)}\in V_1(1, \varepsilon)$. Setting
$$\widetilde{W}_1(u)= \dis - (\sum_{k=1}^n b_k(y_1))Z_{1}(u),$$
where $Z_{1}(u)$ is defined in \eqref{2.A}. Using expansion $(i)$ of Proposition \ref{prop4}, we get
$$<\partial J(u), \widetilde{W}_1(u)> =-(n-2) J(u) \alpha_i^2 \frac{c_1}{K(a_1)}\frac{(\sum_{k=1}^n b_k(y_1))^2}{\lambda_1^{\beta(y_1)}}+ O\Big(|a_1-y_1|^{\beta(y_1)}\Big).$$
Observe that $|a_1-y_1|^{\beta(y_1)}= o\Big(\frac{1}{\lambda_1^{\beta(y_1)}}\Big)$, as $\delta$ small. Thus \begin{equation}\label{01}
 <\partial J(u), \widetilde{W}_1(u)> \leq -\frac{c}{\lambda_1^{\beta(y_1)}}.
\end{equation}
Under $(f)_\beta$-condition, it is easy to verify that for any point $a$ near a
critical point $y$ of $K$, we have
\begin{equation}\label{02}
    \frac{|\nabla K(a)|}{\lambda}= O\Big(\frac{1}{\lambda^{\beta(y)}}\Big), \; \mbox{ if $\lambda|a-y|$ is bounded.}
\end{equation}
In our statement, inequality \eqref{01} implies
$$ <\partial J(u), \widetilde{W}_1(u)> \leq -c\bigg(\frac{1}{\lambda_1^{\beta(y_1)}}+ \frac{|\nabla K(a_1)|}{\lambda_1}\bigg).$$
Observe that the component $\lambda_1(s)$ of the motion $u(s)$ of $\widetilde{W}_1$ increases and goes to $\infty$, only if $y_1\in \Gamma^-$.\\

\n $\bullet$ \underline{Pseudogradient in $V_1(p, \varepsilon), p\geq 2$:}\vskip0.15cm

\n In this region, we organize our construction in to three steps by decomposing $V_1(p, \varepsilon)$ on three subsets.

$$W_1(p, \varepsilon)= \{u=\dis \sum_{i=1}^p \alpha_i \delta_{(a_i, \lambda_i)}\in V_1(p, \varepsilon), s.t., \frac{1}{\beta(y_i)}+ \frac{1}{\beta(y_j)}- \frac{2}{n-2}>0, \forall 1\leq i\neq j \leq p\}.$$
$$W_2(p, \varepsilon)= \{u=\dis \sum_{i=1}^p \alpha_i \delta_{(a_i, \lambda_i)}\in V_1(p, \varepsilon), s.t., \frac{1}{\beta(y_i)}+ \frac{1}{\beta(y_j)}- \frac{2}{n-2}\geq 0, \forall 1\leq i\neq j \leq p\ $$$$ \mbox{ and } \exists i_0\neq j_0 \mbox{ staisfying } \frac{1}{\beta(y_{i_0})}+ \frac{1}{\beta(y_{j_0})}- \frac{2}{n-2}=0\}.$$
$$W_3(p, \varepsilon)= \{u=\dis \sum_{i=1}^p \alpha_i \delta_{(a_i, \lambda_i)}\in V_1(p, \varepsilon), s.t., \exists i_0\neq j_0 \mbox{ staisfying } \frac{1}{\beta(y_{i_0})}+ \frac{1}{\beta(y_{j_0})}- \frac{2}{n-2}<0\}.$$
Next, we denote $\beta_i$ the flatness order of $y_i$ instead $\beta(y_i)$.\\

\n $\bullet$ \textbf{Step 1.  Pseudogradient in $W_1(p, \varepsilon)$:}

\n Let $u=\dis \sum_{i=1}^p \alpha_i \delta_{(a_i, \lambda_i)}\in W_1(p, \varepsilon)$. We introduce the following Lemma.
\begin{lem}\label{lem1} For any $1\leq i\neq j\leq p$, we have
$$\varepsilon_{ij}= o\Big(\frac{1}{\lambda_i^{\beta_i}}\Big)+ o\Big(\frac{1}{\lambda_j^{\beta_j}}\Big), \; \mbox{ as } \varepsilon\rightarrow 0.$$
\end{lem}
\begin{pf}
Let $1\leq i\neq j\leq p$. It is easy to verify that
$$\varepsilon_{ij}\sim \frac{1}{(\lambda_i\lambda_j)^\frac{n-2}{2}} \mbox{ and } \beta_i+\beta_j-\frac{2\beta_i\beta_j}{n-2}>0.$$
Let $0<\gamma< \min\Big(\frac{n-2}{2}, \beta_i+\beta_j-\frac{2\beta_i\beta_j}{n-2}\Big)$.

\n If $\lambda_j^\frac{n-2}{2}\geq \dis \frac{1}{\varepsilon^\gamma}\lambda_i^{\beta_i-\frac{n-2}{2}}$, then $\dis \varepsilon_{ij}\leq \varepsilon^\gamma \frac{1}{\lambda_i^{\beta_i}}.$

\n If $\lambda_j^\frac{n-2}{2}\leq \dis \frac{1}{\varepsilon^\gamma}\lambda_i^{\beta_i-\frac{n-2}{2}}$. This implies that  $\beta_i>\dis \frac{n-2}{2}$,  (if not, $\beta_i- \frac{n-2}{2}\leq 0$, therefore $\lambda_j^\frac{n-2}{2}\leq \dis \frac{1}{\varepsilon^\gamma} < \frac{1}{\varepsilon^\frac{n-2}{2}}$ since $\gamma< \frac{n-2}{2}$. Thus $\lambda_j\leq \frac{1}{\varepsilon}$ which is a contradiction). It follows that

$$\lambda_j^\frac{n-2}{2\beta_i-(n-2)}\leq \frac{1}{\varepsilon^{\frac{2\gamma}{2\beta_i-(n-2)}}}\lambda_i.$$
Hence$$\frac{1}{\lambda_i^\frac{n-2}{2}}\leq \frac{1}{\varepsilon^{\frac{(n-2)\gamma}{2\beta_i-(n-2)}}} \frac{1}{\lambda_j^{\frac{n-2}{2}\frac{n-2}{2\beta_i-(n-2)}}}.$$
Thus,

$$\varepsilon_{ij}\leq \frac{1}{\varepsilon^{\frac{(n-2)\gamma}{2\beta_i-(n-2)}}} \frac{1}{\lambda_j^{\frac{n-2}{2}(1+\frac{n-2}{2\beta_i-(n-2)})}}$$
$$\leq\frac{1}{\varepsilon^{\frac{(n-2)\gamma}{2\beta_i-(n-2)}}} \frac{1}{\lambda_j^{\beta_j+ \frac{n-2}{2\beta_i-(n-2)}(\beta_i+\beta_j-\frac{2\beta_i\beta_j}{n-2})}}$$
$$\leq\frac{\varepsilon^{\frac{n-2}{2\beta_i-(n-2)}(\beta_i+\beta_j-\frac{2\beta_i\beta_j}{n-2})}}
{\varepsilon^{\frac{(n-2)\gamma}{2\beta_i-(n-2)}}} \frac{1}{\lambda_j^{\beta_j}}, \mbox{ since } \frac{1}{\lambda_j}< \varepsilon.$$
Thus,
$$\varepsilon_{ij} \leq\varepsilon^{\frac{n-2}{2\beta_i-(n-2)}(\beta_i+\beta_j-\frac{2\beta_i\beta_j}{n-2}-\gamma)} \frac{1}{\lambda_j^{\beta_j}}.$$
Using the fact tat $\gamma< \beta_i+\beta_j-\frac{2\beta_i\beta_j}{n-2}$, Lemma \ref{lem1} follows.
\end{pf}

\n The construction of the required pseudogradient in $W_1(p, \varepsilon)$ depends to the two following cases.

\n $\bullet$ Case 1, $\forall 1\leq i \leq p, \; - \sum_{k=1}^n b_k(y_i)>0$. Using the notation \eqref{2.A}, we set
$$W_1^1(u)= \sum_{i=1}^p Z_i(u).$$

\n The first expansion of Proposition \ref{prop4} yields

$$ <\partial J(u), {W}_1^1(u)> \leq c \sum_{i=1}^p \frac{(\sum_{k=1}^n b_k(y_i))}{\lambda_i^{\beta_i}}+  \sum_{j\neq i} O\Big(\varepsilon_{ij}\Big).$$
Using the estimate of Lemma \ref{lem1} and the fact that $\sum_{k=1}^n b_k(y_i)<0, \forall i$, we get

$$ <\partial J(u), {W}_1^1(u)> \leq -c \bigg(\sum_{i=1}^p \frac{1}{\lambda_i^{\beta_i}}+  \sum_{j\neq i} \varepsilon_{ij}\bigg)$$$$\leq -c \bigg(\sum_{i=1}^p\bigg( \frac{1}{\lambda_i^{\beta_i}}+\frac{|\nabla K(a_i)|}{\lambda_i}\bigg)+  \sum_{j\neq i} \varepsilon_{ij}\bigg),$$
since by \eqref{02}, $\dis \frac{|\nabla K(a_i)|}{\lambda_i}= O\Big(\frac{1}{\lambda_i^{\beta_i}}\Big), \forall i.$

\n In this region, we observe that a concentration phenomenon happens, since all $\lambda_i, i=1, \ldots, p$, increase and tend to $+\infty$.

\n $\bullet$ Case 2, $\exists   i_1, \ldots,   i_q$ such that $- \sum_{k=1}^n b_k(y_i)<0, \forall j=1, \ldots, q$. Let

$$I= \{i, 1\leq i \leq p, s.t., \lambda_i^{\beta_i}\leq \frac{1}{2} \min_{1\leq j \leq q} \lambda_{i_j}^{\beta_{i_j}}\}.$$
For any $i\in I, (-\sum_{k=1}^n b_k(y_i))>0$. Therefore,
$$ <\partial J(u), \sum_{i\in I}Z_i(u)> \leq -c \sum_{i\in I} \frac{1}{\lambda_i^{\beta_i}}+  O\Big(\sum_{i\in I, j\not \in I} \varepsilon_{ij}\Big)$$
\begin{equation}\label{03}
\leq -c \sum_{i\in I}\Big( \frac{1}{\lambda_i^{\beta_i}}+ \frac{|\nabla K(a_i)|}{\lambda_i}\Big)+  o\Big(\sum_{i\not \in I} \frac{1}{\lambda_i^{\beta_i}}\Big).
\end{equation}
From another part, we decrease all the $\lambda_{i_{j}}, \ j=1,...,q$ with respect to the vector field $-Z_{i_{j}}(u)$ defined in \eqref{2.A}. We have

$$ <\partial J(u), -\sum_{j=1}^q Z_{i_{j}}(u)> \leq -c \sum_{j=1}^q \frac{1}{\lambda_{i_j}^{\beta_{i_j}}}+  O\Big(\sum_{j=1, k\neq i_j}^q \varepsilon_{i_jk}\Big)$$
\begin{equation}\label{04}
\leq -c \sum_{i\not\in I}\Big( \frac{1}{\lambda_i^{\beta_i}}+ \frac{|\nabla K(a_i)|}{\lambda_i}\Big)+  o\Big(\sum_{i\in I} \frac{1}{\lambda_i^{\beta_i}}\Big).
\end{equation}
Let $$W_1^2(u)= \sum_{i\in I}Z_i(u)- \sum_{j=1}^q Z_{i_{j}}(u).$$
$W_1^2$ satisfies $(P.S)$ condition in this region. Moreover, by \eqref{03}, \eqref{04} and Lemma \ref{lem1}, we have
$$<\partial J(u), {W}_1^2(u)> \leq -c \bigg(\sum_{i=1}^p\Big(\frac{1}{\lambda_i^{\beta}}+ \frac{|\nabla K(a_i)|}{\lambda_i}\Big)+ \sum_{j\neq i}\varepsilon_{ij}\bigg).$$
At this step, the required pseudogradient in $W_1(p, \varepsilon)$ say $W_1$ is aconvex combination of $W_1^1$ and $W_1^2$.\\

\n $\bullet$ \textbf{Step 2.  Pseudogradient in $W_2(p, \varepsilon)$:}

\n Let $u= \sum_{i=1}^p \alpha_i \delta_{(a_i, \lambda_i)}\in W_2(p, \varepsilon)$. Setting
$$A_0=\{i, 1 \leq i \leq p, s.t., \exists j, 1\leq j \leq p \mbox{ satisfying } \frac{1}{\beta_i}+\frac{1}{\beta_j}- \frac{2}{n-2}=0 \}.$$
We distinguish two cases.

\n $\bullet$  \underline{Case 1.} $\exists i \in A_0$ such that $\beta_i=n-2$. In this case $\beta_i=n-2, \forall i\in A_0.$

\n We decompose $u$ as follows:

$$u= \sum_{i\not\in A_0} \alpha_i \delta_{(a_i, \lambda_i)}+  \sum_{i\in A_0} \alpha_i \delta_{(a_i, \lambda_i)}= u_1+u_2.$$
The following Lemma  is extracted from (\cite{cbd1}, Proposition 4.3.1).

\begin{lem}\label{lem2}\cite{cbd1} Under assumption $\mathbf{(H_1)}$, there exists a bounded pseudogradient $Y$ in $C_{n-2}(r, \varepsilon):=\{\dis u= \sum_{i=1}^r \alpha_i \delta_{(a_i, \lambda_i)}, \; \lambda_i|a_i-y_i|< \delta, \beta_i=n-2, \forall i=1, \ldots, r, \; y_i\neq y_j, \forall 1\leq i\neq j \leq r\}$, $r\geq 1$ such that for any $u= \sum_{i=1}^r \alpha_i \delta_i \in C_{n-2}(r, \varepsilon)$, we have

$$<\partial J(u), Y(u)> \leq -c \bigg(\sum_{i=1}^r\Big(\frac{1}{\lambda_i^{n-2}}+ \frac{|\nabla K(a_i)|}{\lambda_i}\Big)+ \sum_{j\neq i}\varepsilon_{ij}\bigg).$$
In addition, the only case where the maximum of $\lambda_i$'s is not bounded is when $y_i\in \Gamma^-, \forall i=1, \ldots, r$ and  $\rho(y_1, \ldots, y_r)>0$.
\end{lem}

\n The construction of $W_2^1$; the required pseudogradient in this case, depends on the following three subcases. Let us denote $q= \sharp A_0^c$ and $r= \sharp A_0$.

\n $\bullet$ Subcase 1. $u_1\in W_1^\infty(q, \varepsilon):= \{\dis u=\sum_{i=1}^q \alpha_i \delta_{(a_i, \lambda_i)}\in W_1(q, \varepsilon), s.t., -\sum_{k=1}^n b_k(y_i)>0, \forall i=1, \ldots, q\}$ and $u_2\in C_{n-2}^\infty(r, \varepsilon):= \{u=\dis \sum_{i=1}^r \alpha_i \delta_{(a_i, \lambda_i)}\in C_{n-2}(r, \varepsilon), s.t., y_i\in \Gamma^-, \forall i=1, \ldots, r$ and $\rho(y_1, \ldots, y_r)>0\}$.

\n In this region, $Y(u_2)= \dis \sum_{i\in A_0}\alpha_i \lambda_i \frac{\partial \delta_i}{\partial \lambda_i}$, see \cite{cbd1}. It satisfies

$$<\partial J(u), Y(u_2)> \leq -c \bigg(\sum_{i\in A_0}\Big(\frac{1}{\lambda_i^{\beta_i}}+ \frac{|\nabla K(a_i)|}{\lambda_i}\Big)+ \sum_{i\neq j\in A_0}\varepsilon_{ij}\bigg) + \sum_{i\in A_0, j\not\in A_0} O\Big(\varepsilon_{ij}\Big).$$
From another part, for $W_1^1(u_1)$, where $W_1^1$ is   the vector field defined in $W_1(q, \varepsilon)$, we have

$$<\partial J(u), W_1^1(u_1)> \leq -c \sum_{i\not \in A_0}\Big(\frac{1}{\lambda_i^{\beta_i}}+ \frac{|\nabla K(a_i)|}{\lambda_i}\Big) + \sum_{i\not\in A_0, j\neq i} O\Big(\varepsilon_{ij}\Big).$$
Let $V_1^1(u)= W_1(u_1)+ Y(u_2)$. From the above two inequalities, we get
$$<\partial J(u), V_1^1(u)> \leq -c \bigg(\sum_{i=1}^p\Big(\frac{1}{\lambda_i^{\beta_i}}+ \frac{|\nabla K(a_i)|}{\lambda_i}\Big)+ \sum_{i\neq j}\varepsilon_{ij}\bigg),$$
since by Lemma \ref{lem1}
$$\varepsilon_{ij}= o\Big(\frac{1}{\lambda_i^{\beta_i}}\Big)+ o\Big(\frac{1}{\lambda_j^{\beta_j}}\Big), \mbox{ if } i \mbox{ or } j \not\in A_0.$$
Observe that along the flow  lines of $V_1^1$, $\lambda_i(s)$ moves according the differential equation $\dot{\lambda}_i= \lambda_i, \forall i=1, \ldots, p$. It is therefore a concentration phenomenon statement.

\n $\bullet$ Subcase 2. $u_1\not\in W_1^\infty(q, \varepsilon)$.

\n We apply $W_1^2$ defined in $W_1(q, \varepsilon)$. It satisfies $(P.S)$ condition on its flow lines and

$$<\partial J(u), W_1^2(u_1)> \leq -c \bigg(\sum_{i\not\in A_0} \Big(\frac{1}{\lambda_i^{\beta_i}}+ \frac{|\nabla K(a_i)|}{\lambda_i}\Big)+ \sum_{i\neq j\not \in A_0}\varepsilon_{ij}\bigg)+O\Big(\sum_{i\not\in A_0,  j \in A_0}\varepsilon_{ij}\Big).$$
Let $i_1$ be an index such that

$$\lambda_{i_1}^{\beta_{i_1}}= \min\{\lambda_{i}^{\beta_i}, i\not\in A_0\}.$$Setting
$$I_1=\{i, 1\leq i \leq p, s.t., \lambda_i^{\beta_i}\leq\frac{1}{2}\lambda_{i_1}^{\beta_{i_1}}\}.$$
The above inequality yields

$$<\partial J(u), W_1^2(u_1)> \leq -c \bigg(\sum_{i\not\in I_1} \Big(\frac{1}{\lambda_i^{\beta_i}}+ \frac{|\nabla K(a_i)|}{\lambda_i}\Big)+ \sum_{i\neq j\not \in A_0}\varepsilon_{ij}\bigg)+\sum_{i\in I_1} o\Big(\frac{1}{\lambda_i^{\beta_i}}\Big).$$
In order to make appear $\dis-\sum_{i\not\in I_1, j\neq i}\varepsilon_{ij}$ in the above upper bound, we decrease all $\lambda_i, i\not\in I_1$. Using the fact that

$$\lambda_i \frac{\partial \varepsilon_{ij}}{\partial \lambda_i}\leq - \varepsilon_{ij}, \mbox{ if } |a_i-a_j|\geq d_0>0,$$we obtain

$$<\partial J(u), - \sum_{i\not\in I_1}Z_i(u)> \leq -c \sum_{i\not\in I_1, j\neq i} \varepsilon_{ij}+\sum_{i\not \in I_1} O\Big(\frac{1}{\lambda_i^{\beta_i}}\Big).$$
For a positive constant $m$ small enough, we get

$$<\partial J(u), W_1^2(u_1)-m\sum_{i\not\in I_1}Z_i(u)> \leq -c \bigg(\sum_{i\not\in I_1} \Big(\frac{1}{\lambda_i^{\beta_i}}+ \frac{|\nabla K(a_i)|}{\lambda_i}\Big)+ \sum_{i\not\in I_1, j\neq i}\varepsilon_{ij}\bigg)+\sum_{i\in I_1} o\Big(\frac{1}{\lambda_i^{\beta_i}}\Big).$$
Let $\hat{u}= \dis \sum_{i\in I_1}\alpha_i \delta_{(a_i, \lambda_i)}$. In our statement, we have $\beta_i=n-2$ for any $i\in I_1$. We apply $Y(\hat{u})$ defined in $C_{n-2}(\sharp I_1, \varepsilon)$. It satisfies
$$<\partial J(u), Y(\hat{u})> \leq -c \bigg(\sum_{i\in I_1} \Big(\frac{1}{\lambda_i^{\beta_i}}+ \frac{|\nabla K(a_i)|}{\lambda_i}\Big)+ \sum_{i\neq j\in I_1}\varepsilon_{ij}\bigg)+\sum_{i\in I_1, j\not\in I_1} O\Big(\varepsilon_{ij}\Big).$$
For $m'>0$ and small, let

$$V_1^2(u)=W_1^2(u_1)-m\sum_{i\not\in I_1}Z_i(u)+ m'Y(\hat{u}).$$
$V_1^2$ satisfies $(P.S)$ condition in this region and
$$<\partial J(u), V_1^2(u)> \leq -c \bigg(\sum_{i=1}^p \Big(\frac{1}{\lambda_i^{\beta_i}}+ \frac{|\nabla K(a_i)|}{\lambda_i}\Big)+ \sum_{i\neq j}\varepsilon_{ij}\bigg).$$

\n $\bullet$ Subcase 3. $u_2\not\in C_{n-2}^\infty(r, \varepsilon)$.

\n We apply $Y(u_2)$ defined in $C_{n-2}^\infty(r, \varepsilon)$. It satisfies $(P.S)$ condition and

$$<\partial J(u), Y(u_2)> \leq -c \bigg(\sum_{i\in A_0} \Big(\frac{1}{\lambda_i^{\beta_i}}+ \frac{|\nabla K(a_i)|}{\lambda_i}\Big)+ \sum_{i\neq j\in A_0}\varepsilon_{ij}\bigg)+ \sum_{i\in A_0, j\not\in A_0}O\Big(\varepsilon_{ij}\Big).$$
We proceed as the above subcase, we derive that for $I_1=\dis\{i, 1\leq i\leq p, s.t., \lambda_i^{\beta_i} \leq\frac{1}{2} \min_{i\in A_0} \lambda_i^{\beta_i} \}$, $\hat{u}=\dis \sum_{i\in I_1}\alpha_i \delta_i $ and $V_1^3(u)= Y(u_2)-m \sum_{i\not\in I_1}Z_i(u)+ m' W_1(\hat{u})$,

$$<\partial J(u), V_1^3(u)> \leq -c \bigg(\sum_{i=1}^p \Big(\frac{1}{\lambda_i^{\beta_i}}+ \frac{|\nabla K(a_i)|}{\lambda_i}\Big)+ \sum_{i\neq j}\varepsilon_{ij}\bigg).$$
Here $W_1(\hat{u})$ is the vector field defined in $W_1(\sharp I_1, \varepsilon)$. $V_1^3$ satisfies $(P.S)$ condition. In this case $W_2^1$ is defined by a convex combination of $V_1^1, V_1^2$ and $V_1^3$.\\

\n $\bullet$ Case 2. $\forall i \in A_0, \ \beta_i \neq n-2$. In this case, there exists only one index say $j_0$ belonging to $A_0$ such that $\beta_{j_0}>n-2$. Therefore $A_0$ is reduced to $A_0= \{j_0\}\cup A_{j_0}$, where

$$A_{j_0}=\{i, 1\leq i \leq p, s.t., \beta_i< n-2 \mbox{ and } \frac{1}{\beta_i}+\frac{1}{\beta_{j_0}}-\frac{2}{n-2}=0\}.$$
We introduce the following Lemma

\begin{lem}\label{lem3}
Under assumption $\mathbf{(H_2)}$ there exists a bounded pseudogradient $\widetilde{V}(u)$ satisfying $(P.S)$ condition and
$$<\partial J(u), \widetilde{V}(u)> \leq -c \sum_{i\in A_{j_0}} \Big(\frac{1}{\lambda_i^{\beta_i}}+ \varepsilon_{ij_0}\Big)+ o\Big(\sum_{i\not \in A_{j_0}}\frac{1}{\lambda_i^{\beta_i}}\Big).$$
\end{lem}

\begin{pf}
For any $1\leq i \neq j\leq p$, we have $y_i\neq y_j$, therefore
$$\varepsilon_{ij}= \bigg(\frac{2}{(1-\cos d(a_i, a_j))\lambda_i \lambda_j}\bigg)^{\frac{n-2}{2}}(1+o(1))$$$$=2^\frac{n-2}{2} \frac{G(a_i, a_j)}{(\lambda_i \lambda_j)^\frac{n-2}{2}}\Big(1+o(1)\Big)=2^\frac{n-2}{2} \frac{G(y_i, y_j)}{(\lambda_i \lambda_j)^\frac{n-2}{2}}\Big(1+o(1)\Big).$$
It follows that

$$\lambda_i \frac{\partial \varepsilon_{ij}}{\partial \lambda_i}=-\frac{n-2}{2}2^\frac{n-2}{2} \frac{G(y_i, y_j)}{(\lambda_i \lambda_j)^\frac{n-2}{2}}\Big(1+o(1)\Big).$$
The expansion of $\dis \lambda_i \frac{\partial J(u)}{\partial \lambda_i}$ is then reduced to
$$<\partial J(u), \alpha_i \lambda_i \frac{\partial \delta_i}{\partial \lambda_i}> =(n-2) J(u)\bigg( \tilde{c}_1 \sum_{j\neq i}\alpha_i \alpha_j 2^\frac{n-2}{2} \frac{G(y_i, y_j)}{(\lambda_i \lambda_j)^\frac{n-2}{2}}+ \alpha_i^2 \frac{c_i}{K(y_i)}\frac{\sum_{k=1}^n b_k(y_i)}{\lambda_i^{\beta_i}}\bigg)$$$$+ o\Big(\frac{1}{\lambda_i^{\beta_i}}\Big)+ o\Big(\sum_{j\neq i}\varepsilon_{ij}\Big).$$
Using the fact that $\dis \alpha_i^\frac{4}{n-2} J(u)^\frac{n}{n-2}K(a_i)=1+o(1)$, we get
$$\alpha_i \alpha_j = \frac{1}{\Big( K(a_i)K(a_j)\Big)^\frac{n-2}{4}}J(u)^\frac{-n}{2}(1+o(1))\; \mbox{ and } \alpha_i^2= \frac{1}{K(a_i)^\frac{n-2}{2}}J(u)^\frac{-n}{2}(1+o(1)).$$
Therefore,
$$<\partial J(u), \alpha_i \lambda_i \frac{\partial \delta_i}{\partial \lambda_i}> =(n-2) J(u)^{1-\frac{n}{2}} \bigg( \sum_{j\neq i} \tilde{c}_1  2^{\frac{n-2}{2}} \frac{G(y_i, y_j)}{\Big( K(y_i) K(y_j)\Big)^\frac{n-2}{4}}\frac{1}{(\lambda_i \lambda_j)^\frac{n-2}{2}}$$$$+  \frac{c_i}{K(y_i)^\frac{n}{2}}\frac{\sum_{k=1}^n b_k(y_i)}{\lambda_i^{\beta_i}}\bigg)+ o\Big(\frac{1}{\lambda_i^{\beta_i}}\Big)+ o\Big(\sum_{j\neq i}\varepsilon_{ij}\Big).$$
In order to simplify our notations, we denote in the next,

$$K_i= c_i \frac{\sum_{k=1}^n b_k(y_i)}{K(y_i)^\frac{n}{2}}, i=1, \ldots, p$$
and $$K_{ij}=  \tilde{c}_1 2^\frac{n-2}{2}  \frac{G(y_i, y_j)}{\Big( K(y_i) K(y_j)\Big)^\frac{n-2}{4}}, 1\leq i \neq j\leq p.$$Let $\gamma$ be a positive constant small enough. We distinguish two subcases.\\

\n $\bullet$ Subcase 1. $\exists i_0\in A_{j_0}$ such that
\begin{equation}\label{2.1}
\frac{K_{i_0j_0}}{(\lambda_{i_0}\lambda_{j_0})^\frac{n-2}{2}}< (1-\gamma) \frac{|K_{i_0}|}{\lambda_{i_0}^{\beta_{i_0}}}
\end{equation}
In this region we set
$$V_{i_0}(u)= \Big(-\sum_{k=1}^n b_k(y_{i_0})\Big)Z_{i_0}(u).$$
We have
$$<\partial J(u), V_{i_0}(u)> =(n-2) J(u)^{1-\frac{n}{2}}\bigg( \sum_{j\neq i_0} \frac{(-\sum_{k=1}^n b_k(y_{i_0})) K_{i_0 j}}{(\lambda_{i_0} \lambda_j)^\frac{n-2}{2}}- \frac{\sum_{k=1}^n b_k(y_{i_0})K_{i_0}}{\lambda_{i_0}^{\beta_{i_0}}}\bigg)$$$$+ o\Big(\frac{1}{\lambda_{i_0}^{\beta_{i_0}}}\Big)+ o\Big(\sum_{j\neq i_0}\varepsilon_{i_0j}\Big).$$
Observe that for any $j\neq j_0$, we have $\dis \frac{1}{\beta_j}+ \frac{1}{\beta_{i_0}}-\frac{2}{n-2}>0$. Therefore, by Lemma \ref{lem1}, we have

$$\varepsilon_{i_0 j}\sim \frac{1}{(\lambda_{i_0}\lambda_j)^\frac{n-2}{2}}= o\Big(\frac{1}{\lambda_{i_0}^{\beta_{i_0}}}\Big)+ o\Big(\frac{1}{\lambda_{j}^{\beta_{j}}}\Big), \forall j\neq j_0.$$It follows that

$$<\partial J(u), V_{i_0}(u)> =(n-2) J(u)^{1-\frac{n}{2}}\bigg(  \frac{(-\sum_{k=1}^n b_k(y_{i_0})) K_{i_0 j_{0}}}{(\lambda_{i_0} \lambda_{j_0})^\frac{n-2}{2}}- \frac{\sum_{k=1}^n b_k(y_{i_0})K_{i_0}}{\lambda_{i_0}^{\beta_{i_0}}}\bigg)$$$$+ o\Big(\sum_{i\neq i_0, i\neq j_0}\frac{1}{\lambda_{i}^{\beta_{i}}}\Big).$$
It is easy to see that if $(-\sum_{k=1}^n b_k(y_{i_0}))<0$, then

$$<\partial J(u), V_{i_0}(u)> \leq-c\bigg(  \frac{1}{\lambda_{i_0}^{\beta_{i_0}}} + \varepsilon_{i_0 j_0}\bigg)+ o\Big(\sum_{i\neq i_0, i\neq j_0}\frac{1}{\lambda_{i}^{\beta_{i}}}\Big).$$
Now, if $(-\sum_{k=1}^n b_k(y_{i_0}))>0$, from  \eqref{2.1} we obtain

$$<\partial J(u), V_{i_0}(u)> \leq (n-2) J(u)^{1-\frac{n}{2}}\bigg( (1-\gamma) \frac{(-\sum_{k=1}^n b_k(y_{i_0})) |K_{i_0}|}{\lambda_{i_0}^{\beta_{i_0}}}- \frac{\sum_{k=1}^n b_k(y_{i_0})K_{i_0}}{\lambda_{i_0}^{\beta_{i_0}}}\bigg)$$$$+ o\Big(\sum_{i\neq i_0, i\neq j_0}\frac{1}{\lambda_{i}^{\beta_{i}}}\Big)$$$$ \leq (n-2) J(u)^{1-\frac{n}{2}}\bigg( (1-\gamma) c_{i_0} \frac{(-\sum_{k=1}^n b_k(y_{i_0})) |\sum_{k=1}^n b_k(y_{i_0})|}{K(y_{i_0})^\frac{n}{2}\lambda_{i_0}^{\beta_{i_0}}}- c_{i_0} \frac{(\sum_{k=1}^n b_k(y_{i_0}))^2 }{K(y_{i_0})^\frac{n}{2} \lambda_{i_0}^{\beta_{i_0}}}\bigg)$$$$+ o\Big(\sum_{i\neq i_0, i\neq j_0}\frac{1}{\lambda_{i}^{\beta_{i}}}\Big)$$
$$\leq -\gamma c \frac{(\sum_{k=1}^n b_k(y_{i_0}))^2 }{\lambda_{i_0}^{\beta_{i_0}}}+ o\Big(\sum_{i\neq i_0, i\neq j_0}\frac{1}{\lambda_{i}^{\beta_{i}}}\Big).$$
Using \eqref{2.1}, we get
$$<\partial J(u), V_{i_0}(u)>\leq -c\bigg(  \frac{1}{\lambda_{i_0}^{\beta_{i_0}}} + \varepsilon_{i_0 j_0}\bigg)+ o\Big(\sum_{i\neq i_0, i\neq j_0}\frac{1}{\lambda_{i}^{\beta_{i}}}\Big).$$
Setting

$$I= \{i\in A_0, s.t., \frac{K_{ij_0}}{(\lambda_i \lambda_{j_0})^\frac{n-2}{2}}<(1-\gamma)\frac{|K_i|}{\lambda_i^{\beta_i}}\}.$$
For any $i\in I$, we set $\dis V_i(u)=  (-\sum_{k=1}^n b_k(y_{i}))Z_i(u)$. It follows  from the above computation that

\begin{equation}\label{2.1'}
<\partial J(u), \sum_{i\in I}V_{i}(u)>\leq -c \sum_{i\in I}\bigg(  \frac{1}{\lambda_{i}^{\beta_{i}}} + \varepsilon_{i j_0}\bigg)+ o\Big(\sum_{i\not\in I}\frac{1}{\lambda_{i}^{\beta_{i}}}\Big).
\end{equation}
For $i\in A_{j_0}\setminus I$, we have
\begin{equation}\label{2.1''}
    \frac{K_{ij_0}}{(\lambda_i \lambda_{j_0})^\frac{n-2}{2}}\geq (1-\gamma) \frac{|K_i|}{\lambda_i^{\beta_i}}.
\end{equation}
This implies that
\begin{equation}\label{2.2}
   \lambda_i^{\beta_i-\frac{n-2}{2}}\geq (1-\gamma) \frac{|K_i|}{K_{ij_0}}\lambda_{j_0}^{\frac{n-2}{2}}\geq c \lambda_{j_0}^{\frac{n-2}{2}},
\end{equation}
for some positive constant $c$. Using now the fact that $i_0$ satisfies \eqref{2.1}, we get
\begin{equation}\label{2.3}
   \lambda_{j_0}^{\frac{n-2}{2}}> \frac{1}{1-\gamma} \frac{K_{i_0j_0}}{|K_{i_0}|}\lambda_{i_0}^{\beta_{i_0}-\frac{n-2}{2}}\geq c' \lambda_{i_0}^{\beta_{i_0}-\frac{n-2}{2}}.
\end{equation}
\eqref{2.2} and \eqref{2.3} yield
\begin{equation}\label{2.4}
  \lambda_{i}^{\beta_{i}-\frac{n-2}{2}}\geq c'' \lambda_{i_0}^{\beta_{i_0}-\frac{n-2}{2}}.
\end{equation}
Observe that for any $i\in A_{j_0}$, $\beta_i= \beta_{i_0}>\frac{n-2}{2}$. This with \eqref{2.4} imply that $\dis \frac{1}{\lambda_i^{\beta_i}}\leq M \frac{1}{\lambda_{i_0}^{\beta_{i_0}}}$, for some positive constant $M$. We therefore obtain from \eqref{2.1'}
$$
<\partial J(u), \sum_{i\in I}V_{i}(u)>\leq -c \bigg(\sum_{i\in A_{j_0}} \frac{1}{\lambda_{i}^{\beta_{i}}} + \sum_{i\in I}\varepsilon_{i j_0}\bigg)+ o\Big(\sum_{i\not\in A_{j_0}}\frac{1}{\lambda_{i}^{\beta_{i}}}\Big).
$$

\n In order to make appear $\dis-\sum_{i\in A_{j_0}}\varepsilon_{ij_0}$ in the above upper bound, we decrease all $\lambda_i, i\in A_{j_0}$. Since $|a_i-a_j|\geq d>0$ for any $i \neq j$, therefore  $\lambda_i \dis\frac{\partial \varepsilon_{ij}}{\partial \lambda_i}\leq - \varepsilon_{ij}$. Thus by Proposition \ref{prop4}, we obtain
$$<\partial J(u),- \sum_{i\in A_{j_0}}Z_i(u)>\leq -c \sum_{i\in A_{j_0}, i\neq j} \varepsilon_{i j}+\sum_{i\in A_{j_0}} O\Big(\frac{1}{\lambda_{i}^{\beta_{i}}}\Big).$$
Let $m>0$ small enough and let
$$\widetilde{V}_1= \sum_{i\in I}V_i-m \sum_{i\in A_{j_0}}Z_i.$$
From the above two estimates, we have
$$<\partial J(u), \widetilde{V}_{1}(u)>\leq -c \sum_{i\in A_{j_0}} \bigg(\frac{1}{\lambda_{i}^{\beta_{i}}} + \varepsilon_{i j_0}\bigg)+ o\Big(\sum_{i\not\in A_{j_0}}\frac{1}{\lambda_{i}^{\beta_{i}}}\Big).$$
 By construction $\widetilde{V}_1$ satisfies $(P.S)$ condition on its flow lines in this region, since for any $i\in I, \lambda_i^{\beta_i-\frac{n-2}{2}}\leq M \lambda_{j_0}^\frac{n-2}{2}$ and $\lambda_{j_0}$ does not move.\\

\n $\bullet$ Subcase 2. For any $i\in A_{j_0}$,
$$\frac{K_{ij_0}}{(\lambda_{i}\lambda_{j_0})^\frac{n-2}{2}}> (1-2\gamma) \frac{|K_{i}|}{\lambda_{i}^{\beta_{i}}}.$$
Let
$$I=\{i\in A_{j_0}, (1-2\gamma) \frac{|K_i|}{\lambda_i^{\beta_i}}< \frac{K_{ij_0}}{(\lambda_i\lambda_{j_0})^\frac{n-2}{2}}< (1+2\gamma) \frac{|K_i|}{\lambda_i^{\beta_i}}\}.$$
For any $i\in I$, we have

$$(1-2\gamma)^{\frac{2\beta_{j_0}}{n-2}}\bigg( \frac{|K_i|}{K_{ij_0}}\bigg)^\frac{2\beta_{j_0}}{n-2} \frac{1}{\lambda_i^{\beta_i}}< \frac{1}{\lambda_{j_0}^{\beta_{j_0}}}< (1+2\gamma)^\frac{2\beta_{j_0}}{n-2} \bigg( \frac{|K_i|}{K_{ij_0}}\bigg)^\frac{2\beta_{j_0}}{n-2} \frac{1}{\lambda_i^{\beta_i}},$$
since $\dis \frac{2\beta_i-(n-2)}{n-2}\beta_{j_0}= \beta_i$. Moreover,

\begin{equation}\label{3.1}
   \frac{1}{ (1+2\gamma)^{\frac{2\beta_{j_0}}{n-2}}} \bigg( \frac{K_{ij_0}}{|K_i|}\bigg)^\frac{2\beta_{j_0}}{n-2} \frac{1}{\lambda_{j_0}^{\beta_{j_0}}}< \frac{1}{\lambda_i^{\beta_i}}< \frac{1}{(1-2\gamma)^\frac{2\beta_{j_0}}{n-2}} \bigg( \frac{K_{ij_0}}{|K_i|}\bigg)^\frac{2\beta_{j_0}}{n-2} \frac{1}{\lambda_{j_0}^{\beta_{j_0}}}.
\end{equation}
We now compute $\dis\lambda_{j_0}\frac{\partial J(u)}{\partial \lambda_{j_0}}$.\\
Using the fact that $\dis \frac{1}{\beta_{j_0}}+\frac{1}{\beta_{i}} - \frac{2}{n-2}>0, \forall i\not\in A_{j_0}$, therefore by Lemma \ref{lem1}, $\dis \varepsilon_{ij_0}= o\Big(\frac{1}{\lambda_i^{\beta_i}}\Big)+ o\Big(\frac{1}{\lambda_{j_0}^{\beta_{j_0}}}\Big), \forall i\not\in A_{j_0}$, we obtain,

$$ <\partial J(u), \alpha_{j_0} \lambda_{j_0} \frac{\partial \delta_{j_0}}{\partial \lambda_{j_0}}> = (n-2) J(u)^{1-\frac{n}{2}}\bigg(\sum_{i\in A_{j_0}} \frac{K_{ij_0}}{(\lambda_i\lambda_{j_0})^\frac{n-2}{2}} + \frac{K_{j_0}}{\lambda_{j_0}^{\beta_{j_0}}}\bigg) +  o\Big(\sum_{i\not \in A_{j_0}} \frac{1}{\lambda_{i}^{\beta_{i}}}\Big)$$$$= (n-2) J(u)^{1-\frac{n}{2}}\bigg(\sum_{i\in I} \frac{K_{ij_0}}{(\lambda_i\lambda_{j_0})^\frac{n-2}{2}} + \frac{K_{j_0}}{\lambda_{j_0}^{\beta_{j_0}}}\bigg) +  \sum_{i \in A_{j_0}\setminus I}  O\Big(\varepsilon_{ij_0}\Big) + o\Big(\sum_{i\not \in A_{j_0}} \frac{1}{\lambda_{i}^{\beta_{i}}}\Big).$$
Observe that
$$(1-2\gamma) \sum_{i\in I}\frac{|K_i|}{\lambda_i^{\beta_i}}< \sum_{i\in I}\frac{K_{ij_0}}{(\lambda_i\lambda_{j_0})^\frac{n-2}{2}}< (1+2\gamma)\sum_{i\in I} \frac{|K_i|}{\lambda_i^{\beta_i}}.$$
Using the two inequalities of \eqref{3.1}, we obtain

$$(\ell) \leq \sum_{i\in I}\frac{K_{ij_0}}{(\lambda_i\lambda_{j_0})^\frac{n-2}{2}} \leq (L),$$
where $$(L)= \frac{1+2\gamma}{(1-2\gamma)^{\frac{2\beta_{j_0}}{n-2}}} \sum_{i\in I}|K_i| \bigg(\frac{K_{ij_0}}{|K_i|}\bigg)^\frac{2\beta_{j_0}}{n-2} \frac{1}{\lambda_{j_0}^{\beta_{j_0}}}$$
and $$(\ell)= \frac{1-2\gamma}{(1+2\gamma)^{\frac{2\beta_{j_0}}{n-2}}} \sum_{i\in I}|K_i| \bigg(\frac{K_{ij_0}}{|K_i|}\bigg)^\frac{2\beta_{j_0}}{n-2} \frac{1}{\lambda_{j_0}^{\beta_{j_0}}}.$$
Therefore,
\begin{equation}\label{3.2}
    (\widetilde{\ell})\leq <\partial J(u), \alpha_{j_0} \lambda_{j_0} \frac{\partial \delta_{j_0}}{\partial \lambda_{j_0}}> \leq (\widetilde{L}),
\end{equation}
where $$(\widetilde{L})=  (n-2) J(u)^{1-\frac{n}{2}}\bigg(\frac{1+2\gamma}{(1-2\gamma)^{\frac{2\beta_{j_0}}{n-2}}} \sum_{i\in I}|K_i| \bigg(\frac{K_{ij_0}}{|K_i|}\bigg)^\frac{2\beta_{j_0}}{n-2} + K_{j_0}\bigg)\frac{1}{\lambda_{j_0}^{\beta_{j_0}}}$$$$+ O\Big(\sum_{i \in A_{j_0}\setminus I}  \varepsilon_{ij_0}\Big) + o\Big(\sum_{i\not \in A_{j_0}} \frac{1}{\lambda_{i}^{\beta_{i}}}\Big)$$
and $$(\widetilde{\ell})=  (n-2) J(u)^{1-\frac{n}{2}}\bigg(\frac{1-2\gamma}{(1+2\gamma)^{\frac{2\beta_{j_0}}{n-2}}} \sum_{i\in I}|K_i| \bigg(\frac{K_{ij_0}}{|K_i|}\bigg)^\frac{2\beta_{j_0}}{n-2} + K_{j_0}\bigg)\frac{1}{\lambda_{j_0}^{\beta_{j_0}}} $$$$+ O\Big(\sum_{i \in A_{j_0}\setminus I}  \varepsilon_{ij_0}\Big) + o\Big(\sum_{i\not \in A_{j_0}} \frac{1}{\lambda_{i}^{\beta_{i}}}\Big).$$
Let $S_I= \sum_{i\in I} |K_i| \bigg(\frac{K_{ij_0}}{|K_i|}\bigg)^\frac{2\beta_{j_0}}{n-2}$. By taking $\gamma$ small enough, the sign of $\dis \frac{1-2\gamma}{(1+2\gamma)^{\frac{2\beta_{j_0}}{n-2}}}S_I +K_{j_0}$ is the sign of $\dis  \frac{1+2\gamma}{(1-2\gamma)^{\frac{2\beta_{j_0}}{n-2}}}S_I +K_{j_0}$. It is the sign of $S_I+ K_{j_0}$. Recall by assumption $\mathbf{(H_2)}$, we have $S_I+K_{j_0}\neq 0$. Let

$$X(u)= -\mbox{ sign } (S_I+ K_{j_0})\alpha_{j_0}\lambda_{j_0}\frac{\partial \delta_{j_0}}{\partial \lambda_{j_0}}.$$Using \eqref{3.2} we obtain
$$< \partial J(u), X(u)> \leq -\frac{c}{\lambda_{j_0}^{\beta_{j_0}}}+ \sum_{i\in A_{j_0}\setminus I}O(\varepsilon_{ij_0}) + o\Big(\sum_{i\not \in A_{j_0}} \frac{1}{\lambda_{i}^{\beta_{i}}}\Big).$$This with \eqref{3.1} yield,
$$< \partial J(u), X(u)> \leq -c\sum_{i\in I}\frac{1}{\lambda_{i}^{\beta_i}}+ \sum_{i\in A_{j_0}\setminus I}O(\varepsilon_{ij_0}) + o\Big(\sum_{i\not \in A_{j_0}} \frac{1}{\lambda_{i}^{\beta_{i}}}\Big).$$
Of course for any $i\in I$, $ -\frac{1}{\lambda_{i}^{\beta_i}}\leq -c \varepsilon_{ij_0}$, therefore,
\begin{equation}\label{3.3}
< \partial J(u), X(u)> \leq -c\sum_{i\in I}\bigg(\frac{1}{\lambda_{i}^{\beta_i}}+ \varepsilon_{ij_0}\bigg) + \sum_{i\in A_{j_0}\setminus I}O(\varepsilon_{ij_0}) + o\Big(\sum_{i\not \in A_{j_0}} \frac{1}{\lambda_{i}^{\beta_{i}}}\Big).
\end{equation}
$X$ satisfies $(P.S)$ condition in this region, since by \eqref{3.1}, we have $\dis \lambda_{j_0}^{\beta_{j_0}}\leq M \lambda_{i}^{\beta_{i}}$, for some constant $M$ and $\lambda_i$ does not move.

\n For any $i\in A_{j_0}\setminus I$, we have
\begin{equation}\label{3.4}
    \frac{K_{ij_0}}{(\lambda_i \lambda_{j_0})^\frac{n-2}{2}}> (1+2\gamma) \frac{|K_i|}{\lambda_i^{\beta_i}}.
\end{equation}
Therefore,

$$ <\partial J(u), -\alpha_{i} \lambda_{i} \frac{\partial \delta_{i}}{\partial \lambda_{i}}> = (n-2) J(u)^{1-\frac{n}{2}}\bigg(-\sum_{j\neq i} \frac{K_{ij}}{(\lambda_i\lambda_{j})^\frac{n-2}{2}} - \frac{K_{i}}{\lambda_{i}^{\beta_{i}}}\bigg) +  o\Big( \frac{1}{\lambda_{i}^{\beta_{i}}}\Big)+ o\Big(\sum_{j \neq i} \varepsilon_{ij}\Big).$$
For any $j\neq j_0$, $\dis \varepsilon_{ij}= o\Big(\frac{1}{\lambda_{i}^{\beta_{i}}}\Big)  + o\Big(\frac{1}{\lambda_{j}^{\beta_{j}}}\Big).$ Thus,

$$ <\partial J(u), -\alpha_{i} \lambda_{i} \frac{\partial \delta_{i}}{\partial \lambda_{i}}> = (n-2) J(u)^{1-\frac{n}{2}}\bigg(-\frac{K_{ij_0}}{(\lambda_i\lambda_{j_0})^\frac{n-2}{2}} - \frac{K_{i}}{\lambda_{i}^{\beta_{i}}}\bigg) +  o\Big( \sum_{j\neq j_0}\frac{1}{\lambda_{j}^{\beta_{j}}}\Big)$$$$
= (n-2) J(u)^{1-\frac{n}{2}}\bigg(-\frac{2\gamma}{1+2\gamma}\frac{K_{ij_0}}{(\lambda_i\lambda_{j_0})^\frac{n-2}{2}} - \frac{1}{1+2\gamma}\frac{K_{ij_0}}{(\lambda_i\lambda_{j_0})^\frac{n-2}{2}} - \frac{K_{i}}{\lambda_{i}^{\beta_{i}}}\bigg) +  o\Big( \sum_{j\neq j_0}\frac{1}{\lambda_{j}^{\beta_{j}}}\Big).$$
Using \eqref{3.4}, we obtain

$$ <\partial J(u), -\alpha_{i} \lambda_{i} \frac{\partial \delta_{i}}{\partial \lambda_{i}}>\leq -c \frac{2\gamma}{1+2\gamma}\varepsilon_{ij_0} +  o\Big( \sum_{j\neq j_0}\frac{1}{\lambda_{j}^{\beta_{j}}}\Big).$$
By \eqref{3.4} again, we have
$$ <\partial J(u), -\alpha_{i} \lambda_{i} \frac{\partial \delta_{i}}{\partial \lambda_{i}}>\leq -c\Big( \varepsilon_{ij_0}+\frac{1}{\lambda_{i}^{\beta_{i}}}\Big) +  o\Big( \sum_{j\neq j_0}\frac{1}{\lambda_{j}^{\beta_{j}}}\Big).$$
Therefore,
\begin{equation}\label{3.5}
<\partial J(u), - \sum_{i\in A_{j_0}\setminus I}\alpha_{i} \lambda_{i} \frac{\partial \delta_{i}}{\partial \lambda_{i}}>\leq -c \sum_{i\in A_{j_0}\setminus I}\bigg( \frac{1}{\lambda_{i}^{\beta_{i}}}+ \varepsilon_{ij_0}\Big) +  o\Big( \sum_{j\not\in A_{j_0}}\frac{1}{\lambda_{j}^{\beta_{j}}}\Big).
\end{equation}
Let $m$ be a positive constant small enough and let $$\widetilde{V}_2(u) = -\sum_{i\in A_{j_0}\setminus I}\alpha_{i} \lambda_{i} \frac{\partial \delta_{i}}{\partial \lambda_{i}}+ m X(u).$$
\eqref{3.3} and \eqref{3.5} yield
$$
<\partial J(u), \widetilde{V}_2(u)>\leq -c \sum_{i\in A_{j_0}}\bigg( \frac{1}{\lambda_{j}^{\beta_{j}}}+ \varepsilon_{ij_0}\Big) +  o\Big( \sum_{j\not\in A_{j_0}}\frac{1}{\lambda_{j}^{\beta_{j}}}\Big).$$
The required pseudogradient $\widetilde{V}$ of Lemma \ref{lem3} is a convex combination of $\widetilde{V}_1$ and $\widetilde{V}_2$. This finishes the proof of Lemma \ref{lem3}.
\end{pf}

\n We continue our construction of the pseudogradient $W_2^2(u)$ under the assumption that $\beta_i \neq n-2, \ \forall i\in A_0$. Denoting $i_1$ an index of $A_{j_0}$ such that $\dis \lambda_{i_1}= \min \{\lambda_i, i\in A_{j_0}\}$ and setting

$$I_1=\{i, 1\leq i \leq p, s.t., \lambda_i^{\beta_i}\geq \Big(\frac{\lambda_{i_1}}{2}\Big)^{\beta_{i_1}}\}.$$
Of course $A_{j_0}\subset I_1$, since $\beta_i =\beta_{i_1}, \forall i\in A_{j_0}$ and the inequality of Lemma \ref{lem3} becomes

\begin{equation}\label{3.6}
 <\partial J(u), \widetilde{V}(u)>\leq -c \bigg( \sum_{i\in I_{1}}\bigg( \frac{1}{\lambda_{i}^{\beta_{i}}}+ \sum_{i\in A_{j_0}}\varepsilon_{ij_0}\bigg) +  o\Big( \sum_{i\not\in I_{1}}\frac{1}{\lambda_{i}^{\beta_{i}}}\Big).
\end{equation}
Let $\hat{u}= \sum_{i\not\in I_1}\alpha_i \delta_i$. Observe that for any $i, j \in I_1^c$, we have $\dis \frac{1}{\lambda_{i}^{\beta_{i}}}+ \frac{1}{\lambda_{j}^{\beta_{j}}}-\frac{2}{n-2}> 0$. Therefore, $\hat{u}\in W_1(\sharp I_1^c, \varepsilon)$. By applying $W_1$ the pseudogradient defined in $W_1(\sharp I_1^c, \varepsilon)$, we obtain
\begin{equation}\label{3.7}
 <\partial J(u), W_1(\hat{u})>\leq -c \bigg( \sum_{i\not \in I_{1}}\Big( \frac{1}{\lambda_{i}^{\beta_{i}}}+\frac{|\nabla K(a_i)|}{\lambda_i}\Big)+ \sum_{i\neq j \not\in I_{1}}\varepsilon_{ij}\bigg) +  \sum_{i\not\in I_{1}, j\in I_1} O\Big( \varepsilon_{ij}\Big).
\end{equation}
For any $i\not\in A_{j_0}$ and for any $j\neq i$, we have $\dis \frac{1}{\beta_{i}}+ \frac{1}{\beta_{j}} -\frac{2}{n-2}> 0$. Therefore, by Lemma \ref{lem1}, $\dis \varepsilon_{ij}= o\Big(\frac{1}{\lambda_i^{\beta_i}}\Big)+ o\Big(\frac{1}{\lambda_j^{\beta_j}}\Big)$.

\n Let $W_2^2(u)= \dis \widetilde{V}(u) + W_1(\hat{u})$. By \eqref{3.6} and \eqref{3.7}, we have
$$<\partial J(u), W_2^2({u})>\leq -c \bigg( \sum_{i=1}^p\Big( \frac{1}{\lambda_{i}^{\beta_{i}}}+\frac{|\nabla K(a_i)|}{\lambda_i}\Big)+ \sum_{i\neq j}\varepsilon_{ij}\bigg).$$
At this step, the required pseudogradient $W_2$ in $W_2(p, \varepsilon)$ is a convex combination of $ W_2^1$ and $W_2^2$.\\

\n $\bullet $ \textbf{Step 3. Pseudogradient in $W_3(p, \varepsilon)$:}

\n Let $u= \dis\sum_{i=1}^p \alpha_i \delta_i \in W_3(p, \varepsilon)$. We order all the concentration $\lambda_{i}^{\beta_i}$. Without loss of generality, we assume that
$$\lambda_{1}^{\beta_1} \leq \ldots  \leq \lambda_{p}^{\beta_p}.$$Let

$$A_0 = \{i, 1\leq  i \leq p, s.t., \exists j, 1\leq j \leq p, \mbox{ satisfying } \dis \frac{1}{\beta_{i}}+ \frac{1}{\beta_{j}} -\frac{2}{n-2}< 0\}.$$
For $M>0$ and large, we set
$$I = \{i, 1\leq  i \leq p, s.t., \lambda_{i}^{\beta_i} \leq M \lambda_{1}^{\beta_1}\}.$$
We distinguish  two cases.

\n $\bullet $ Case 1. If $\sharp I=1$. In this case we have $ \lambda_{i}^{\beta_i} \geq M \lambda_{1}^{\beta_1}$ for any $i\geq 2$ and thus $ \dis\frac{1}{\lambda_{i}^{\beta_i}}= \Big(\frac{1}{\lambda_{1}^{\beta_1}}\Big)$, as $M$ large. We decrease  all $\lambda_i, i=2, \ldots, p$. Using the first expansion of Proposition \ref{prop4}, we have

$$<\partial J(u), - \sum_{i=2}^p Z_i(u)>\leq -c \sum_{i \neq j} \varepsilon_{ij} + \sum_{i=2}^p O\Big(\frac{1}{\lambda_{i}^{\beta_{i}}}\Big)$$$$\leq -c \sum_{i \neq j} \varepsilon_{ij} + o\Big(\frac{1}{\lambda_{1}^{\beta_{1}}}\Big).$$
Setting $\dot{\lambda}_1= \dis (-\sum_{ k=1}^n b_k(y_1))\lambda_1$. We have

$$<\partial J(u), (-\sum_{ k=1}^n b_k(y_1)) Z_1(u)>\leq -\frac{c}{\lambda_{1}^{\beta_{1}}} + \sum_{j \neq 1} O\Big(\varepsilon_{1j}\Big).$$
For $m> 0$ and small, let\\ $\dis W_3^1(u)= - \sum_{i=2}^p Z_i(u) - m(\sum_{ k=1}^n b_k(y_1))Z_1(u)$. From the above two inequalities,

$$<\partial J(u), W_3^1(u)>\leq - c \bigg(\sum_{j \neq 1} \varepsilon_{1j} + \frac{1}{\lambda_{1}^{\beta_{1}}}\bigg)$$$$\leq - c \bigg(\sum_{i=1}^p \Big(\frac{1}{\lambda_{i}^{\beta_{i}}}+ \frac{|\nabla K(a_i)|}{\lambda_{i}}\Big)+ \sum_{j \neq i} \varepsilon_{1j}\bigg).$$

\n $\bullet $ Case 2. $\sharp I\geq 2$. Our construction depends to the following two subcases.

\n $\bullet $ Subcase 1. $\dis \forall i\neq j \in I, \, \frac{1}{\beta_{i}}+ \frac{1}{\beta_{j}}-\frac{2}{n-2}\geq 0.$

\n We decompose $u$ as follows.

$$u= \hat{u} + \acute{u}= \sum_{i\in I} \alpha_i \delta_i + \sum_{i\not \in I} \alpha_i \delta_i.$$
Observe that $\hat{u}\in W_2(\sharp I, \varepsilon)$. We apply the related pseudogradient $W_2$. We have
\begin{equation}\label{3.8}
  <\partial J(u), W_2(\hat{u})>\leq - c \bigg(\sum_{i\in I} \Big(\frac{1}{\lambda_{i}^{\beta_{i}}}+ \frac{|\nabla K(a_i)|}{\lambda_{i}}\Big)+ \sum_{j \neq i \in I} \varepsilon_{ij}\bigg) + O\Big(\sum_{i \in I, j\not\in I} \varepsilon_{ij}\Big).
\end{equation}
In this region, $W_2$ satisfies $(P.S)$ condition. Since there exists at least an index $i_0 \in A_0$ satisfying $\dis \lambda_{i}^{\beta_{i}}\leq \lambda_{i_0}^{\beta_{i_0}}$ for any $i\in I$. In addition, we have

$$<\partial J(u), - \sum_{i\not\in I}Z_i(u)>\leq -c \sum_{i\not \in I, j\neq i} \varepsilon_{ij} + O\Big(\sum_{ i\not\in I }\frac{1}{\lambda_{i}^{\beta_{i}}}\Big)$$
\begin{equation}\label{3.9}
\leq -c  \sum_{i\not \in I, j\neq i} \varepsilon_{ij} + o\Big(\frac{1}{\lambda_{1}^{\beta_{1}}}\Big).
\end{equation}
Let $\dis W_3^2(u) = mW_2(\hat{u}) - \sum_{i\not\in I}Z_i(u)$, where $m> 0$ and small. We have by \eqref{3.8} and \eqref{3.9}
$$<\partial J(u), W_3^2(u)>\leq  - c \bigg(\sum_{i=1}^p \Big(\frac{1}{\lambda_{i}^{\beta_{i}}}+ \frac{|\nabla K(a_i)|}{\lambda_{i}}\Big)+ \sum_{j \neq i} \varepsilon_{1j}\bigg).$$

\n $\bullet $ Subcase 2. $\dis \exists i_0\neq j_0 \in I, \, \frac{1}{\beta_{i_0}}+ \frac{1}{\beta_{j_0}}-\frac{2}{n-2}< 0.$

\n We introduce the following Lemma.

\begin{lem}\label{lem4}
In our statement, we have

$$\frac{1}{\lambda_{i_0}^{\beta_{i_0}}} = o(\varepsilon_{i_0j_0}), \; \mbox{ as } \varepsilon\rightarrow 0.$$
\end{lem}

\begin{pf}
Let $ \dis 0< \gamma< -\frac{n-2}{2 \beta_{j_0}}\Big(\beta_{i_0} + \beta_{j_0} - \frac{2 \beta_{i_0}\beta_{j_0}}{n-2}\Big).$ We claim that
\begin{equation}\label{3.10}
  \lambda_{j_0}^{\frac{n-2}{2}}< \varepsilon^{\gamma}\lambda_{i_0}^{\beta_{i_0}-\frac{n-2}{2}}.
\end{equation}
Indeed, we argue by contradiction. Assume that

\begin{equation}\label{3.11}
  \lambda_{j_0}^{\frac{n-2}{2}}\geq \varepsilon^{\gamma}\lambda_{i_0}^{\beta_{i_0}-\frac{n-2}{2}}.
\end{equation}
Using the fact that $\dis \beta_{i_0}> \frac{n-2}{2}$ ( as necessary condition for $\dis  \frac{1}{\beta_{i_0}}+ \frac{1}{\beta_{j_0}}-\frac{2}{n-2}< 0$), inequality \eqref{3.11} yields
$$ \lambda_{j_0}^{\beta_{j_0}}\geq \varepsilon^{\frac{2\gamma \beta_{j_0}}{n-2}}\lambda_{i_0}^{\frac{2\beta_{i_0}-(n-2)}{n-2}\beta_{j_0}}$$$$
\geq\varepsilon^{\frac{2\gamma \beta_{j_0}}{n-2}}\lambda_{i_0}^{\beta_{i_0}}\lambda_{i_0}^{-(\beta_{i_0}+\beta_{j_0}-\frac{2\beta_{i_0}\beta_{j_0}}{n-2})}.$$
Using the fact that $\dis \lambda_{i_0}> \frac{1}{\varepsilon}$, we get

$$\lambda_{j_0}^{\beta_{j_0}}\geq \varepsilon^{\gamma`}\lambda_{i_0}^{\beta_{i_0}},$$
where $\gamma`= \dis \frac{2\gamma \beta_{j_0}}{n-2} + (\beta_{i_0}+\beta_{j_0}-\frac{2\beta_{i_0}\beta_{j_0}}{n-2})< 0$. It follows that $\dis \lambda_{i_0}^{\beta_{i_0}} = o(\lambda_{j_0}^{\beta_{j_0}})$ as $\varepsilon\rightarrow 0$, which is a contradiction, since $i_0, j_0 \in I.$ Thus, \eqref{3.10} is valid. Consequently

$$\varepsilon_{i_0j_0}\sim \frac{1}{(\lambda_{i_0}\lambda_{j_0})^\frac{n-2}{2}}> \frac{1}{\varepsilon^{\gamma}}\frac{1}{\lambda_{i_0}^{\beta_{i_0}}}.$$
This conclude the proof of Lemma \ref{lem4}.
\end{pf}

\n By the first expansion of Proposition \ref{prop4}, we have

$$<\partial J(u), - \sum_{i\in I}\alpha_{i} Z_i(u)>\leq -c \sum_{i \in I, j\neq i} \varepsilon_{ij} + O\Big(\sum_{ i\in I }\frac{1}{\lambda_{i}^{\beta_{i}}}\Big).$$
For any $i\in I, \dis \lambda_{i}^{\beta_{i}}\sim \lambda_{i_0}^{\beta_{i_0}}$. Thus by Lemma \ref{lem4},

$$<\partial J(u), - \sum_{i\in I}Z_i(u)>\leq -c \sum_{ i\in I}\bigg( \frac{1}{\lambda_{i}^{\beta_{i}}} + \sum_{ j\neq i} \varepsilon_{ij}\bigg).$$
From another part,

$$<\partial J(u), - \sum_{i\not\in I}Z_i(u)>\leq -c \sum_{i \not \in I, j\neq i} \varepsilon_{ij} + O\Big(\sum_{ i\not \in I }\frac{1}{\lambda_{i}^{\beta_{i}}}\Big)$$$$\leq -c \sum_{i \not \in I, j\neq i} \varepsilon_{ij} + o\Big(\frac{1}{\lambda_{1}^{\beta_{1}}}\Big.)$$
Let $W_3^3(u)= -\dis \sum_{i=1}^p Z_i(u)$. It satisfies $(P.S)$ condition  and
$$<\partial J(u), W_3^3(u)>\leq  - c \bigg(\sum_{i=1}^p \Big(\frac{1}{\lambda_{i}^{\beta_{i}}}+ \frac{|\nabla K(a_i)|}{\lambda_{i}}\Big)+ \sum_{j \neq i} \varepsilon_{1j}\bigg).$$
At this step, the required pseudogradient $W_3$ in $W_3(p, \varepsilon)$ is a convex combination of $W_3^1, W_3^2$ and $W_3^3$.

\n This finishes the construction of $\widetilde{W}_1$ in $V_1(p, \varepsilon)$ which is defined as a convex combination of $W_1, W_2$ and $W_3$.\\

\n $\bullet$ \underline{Pseudogradient in $V_2(p, \varepsilon), p\geq 1$:}
\vskip0.15cm

\n Let $\dis u=\sum_{i=1}^{p}\alpha_i \delta_{(a_i, \lambda_i)}\in V_2(p, \varepsilon)$. Denote $$A_0=\{i, 1\leq i\leq p, s.t., \lambda_i |a_i-y_i|> \delta/2\}.$$We claim the following.

\n $(C_1)$: \; For any $i\in A_0$, there exists a bounded pseudogradient $\widetilde{Y}_i(u)$ satisfying $(P.S)$ condition and $$<\partial J(u), \widetilde{Y}_i(u)> \leq -c \bigg(\frac{1}{\lambda_i^{\beta_i}}+\frac{|\nabla K(a_i)|}{\lambda_i}+ \sum_{j\neq i} \varepsilon_{ij}\bigg).$$
Indeed, let $i\in A_0$.

\n If $\dis \lambda_i |a_i-y_i|\leq \frac{1}{\delta}$, we set
$$X_i(u)= \alpha_i \sum_{k=1}^n b_k \int_{\mathbb{R}^n}\frac{|x_k+\lambda_i(a_i-y_i)_k|^{\beta_i}}{(1+|x|^2)^{n+1}}dx \frac{1}{\lambda_i}\frac{\partial \delta_{(a_i, \lambda_i)}}{\partial (a_i)_k}.$$
$X_i$ is bounded and satisfies $(P.S)$ condition. Using expansion $(ii)'$ of Proposition \ref{prop4}, we have
$$<\partial J(u), X_i(u)> \leq - \frac{c}{\lambda_i^{\beta_i}} \sum_{k=1}^n \bigg(\int_{\mathbb{R}^n}\frac{|x_k+\lambda_i(a_i-y_i)_k|^{\beta_i} x_k}{(1+|x|^2)^{n+1}}dx\bigg)^2
+O\Big(\sum_{j\neq i} \frac{1}{\lambda_i} \Big|\frac{\partial \varepsilon_{ij}}{\partial a_i}\Big|\Big).$$
Using the fact that
$$\sum_{k=1}^n \bigg(\int_{\mathbb{R}^n}\frac{|x_k+\lambda_i(a_i-y_i)_k|^{\beta_i} x_k}{(1+|x|^2)^{n+1}}dx\bigg)\geq c_\delta>0, \mbox{ since } \lambda_i|a_i-y_i|\geq\delta/2$$
and $$\frac{1}{\lambda_i} \Big|\frac{\partial \varepsilon_{ij}}{\partial a_i}\Big|= o(\varepsilon_{ij}), \forall i \neq j, \mbox{ since } |a_i-a_j|\geq d_0,$$
we get
$$<\partial J(u), X_i(u)> \leq -c \bigg(\frac{1}{\lambda_i^{\beta_i}}+\frac{|\nabla K(a_i)|}{\lambda_i}\bigg) + \sum_{j\neq i} o(\varepsilon_{ij}),$$
since by \eqref{02}, we have $\dis \frac{|\nabla K(a_i)|}{\lambda_i}= O\Big(\frac{1}{\lambda_i^{\beta_i}}\Big).$

\n In this statement, let $\dis \widetilde{Y}_i(u)= X_i(u)-m Z_i(u)$, where $m>0$ and small. The preceding inequality and the expansion $(ii)$ of Proposition \ref{prop4}, yield
 $$<\partial J(u), \widetilde{Y}_i(u)> \leq -c \bigg(\frac{1}{\lambda_i^{\beta_i}}+\frac{|\nabla K(a_i)|}{\lambda_i}+ \sum_{j\neq i} \varepsilon_{ij}\bigg)$$
and claim $(C_1)$ valid in this case.

\n If $\dis \lambda_i |a_i-y_i|\geq \frac{1}{\delta}$, we use $Y_i(u)$ the vector field defined in \eqref{2.B}, where the cut-off function $\varphi_{M}$ equals to $\varphi_{\frac{1}{\delta}}$.\\
By $(ii)'$ of Proposition \ref{prop4}, we have

$$< \partial J(u), Y_i(u)>\leq -c \sum_{k=1}^n b_k^2 \frac{|(a_i-y_i)_k|^{\beta_i-1}}{\lambda_i}+ \sum_{j=2}^{[\beta]}O\Big(\frac{|a_i-y_i|^{\beta-j}}{\lambda_i^j}\Big)$$$$+O\Big(\frac{1}{\lambda_i^{\beta_i}}
\Big)+ o(\sum_{j \neq i}\varepsilon_{ij}).$$
For any $j=2, \ldots, [\beta]$, we have
$$\frac{|a_i-y_i|^{\beta-j}}{\lambda_i^{\beta_i}}= o\bigg(\frac{|a_i-y_i|^{\beta-1}}{\lambda_i}\bigg), \mbox{ as $\delta$ small }$$and
$$\frac{1}{\lambda_i^{\beta_i}} = o\bigg(\frac{|a_i-y_i|^{\beta-j}}{\lambda_i}\bigg), \mbox{ as $\delta$ small. }$$
Thus from the preceding inequality, we obtain
$$<\partial J(u), Y_i(u)> \leq -c \bigg(\frac{1}{\lambda_i^{\beta_i}}+\frac{|\nabla K(a_i)|}{\lambda_i}\bigg) + \sum_{j\neq i} o(\varepsilon_{ij}),$$
since $\dis |\nabla K(a_i)|\sim |a_i-y_i|^{\beta-1}$. Let in this case $\dis \widetilde{Y}_i(u)= Y_i(u)-m Z_i(u)$. From the above estimates and expansion $(ii)$ of Proposition \ref{prop4}, we have
$$<\partial J(u), \widetilde{Y}_i(u)> \leq -c \bigg(\frac{1}{\lambda_i^{\beta_i}}+\frac{|\nabla K(a_i)|}{\lambda_i}+ \sum_{j\neq i} \varepsilon_{ij}\bigg).$$
Hence  claim $(C_1)$ follows. Let
 $$\widetilde{A}_0= \dis \{i, 1\leq i \leq p, s.t., \lambda_i^{\beta_i}\geq \frac{1}{2} \min_{j\in A_0} \lambda_j^{\beta_j} \}.$$
\n Using claim $(C_1)$ and the fact that $\dis \frac{|\nabla K(a_i)|}{\lambda_i}= O\Big(\frac{1}{\lambda_i^{\beta_i}}\Big)$ for $i\in \widetilde{A}_0\setminus A_0$, we have

$$<\partial J(u), \sum_{i\in A_0} \widetilde{Y}_i(u)> \leq -c \bigg(\sum_{i\in \widetilde{A}_0}\Big(\frac{1}{\lambda_i^{\beta_i}}+\frac{|\nabla K(a_i)|}{\lambda_i}\Big)+ \sum_{i\in A_0, j\neq i} \varepsilon_{ij}\bigg).$$
Therefore,

$$<\partial J(u), \sum_{i\in A_0} \widetilde{Y}_i(u) - m\sum_{j\in \widetilde{A}_0\setminus A_0}Z_i(u)> \leq -c \sum_{i\in \widetilde{A}_0}\bigg(\frac{1}{\lambda_i^{\beta_i}}+\frac{|\nabla K(a_i)|}{\lambda_i} + \sum_{ j\neq i} \varepsilon_{ij}\bigg).$$
Let  $\dis \hat{u}= \sum_{i\in \widetilde{A}_0^c}\alpha_i \delta_i$ and let $\widetilde{W}_1(\hat{u})$ be the pseudogradient  defined in $V_1(\sharp \widetilde{A_0}^c, \varepsilon)$. We have

$$<\partial J(u), \widetilde{W}_1(\hat{u})> \leq -c \bigg( \sum_{i\in \widetilde{A}_0^c}\bigg(\frac{1}{\lambda_i^{\beta_i}}+\frac{|\nabla K(a_i)|}{\lambda_i}\bigg) + \sum_{ i\neq j \in \widetilde{A}_0^c} \varepsilon_{ij}\bigg) + O\Big( \sum_{ i \in \widetilde{A}_0^c, j\neq i} \varepsilon_{ij}\Big).$$
For $m'>0$ and small, let
$$\widetilde{W}_2(u)= \dis\sum_{i\in A_0} \widetilde{Y}_i(u) +m'\widetilde{W}_1(\hat{u})-m \sum_{j\in \widetilde{A}_0\setminus A_0} Z_i(u).$$
$\widetilde{W}_2$ satisfies $(P.S)$ condition in this region. Moreover it satisfies

$$<\partial J(u), \widetilde{W}_2(u)> \leq -c \bigg(\sum_{i=1}^p \Big(\frac{1}{\lambda_i^{\beta_i}}+\frac{|\nabla K(a_i)|}{\lambda_i}\Big) + \sum_{ j\neq i} \varepsilon_{ij}\bigg).$$

\n $\bullet$ \underline{Pseudogradient in $V_3(p, \varepsilon), p\geq 1$:}
\vskip0.15cm

\n Let $\dis u=\sum_{i=1}^{p}\alpha_i \delta_{(a_i, \lambda_i)}\in V_3(p, \varepsilon)$. For any critical point $y_k$ of $K$, we denote
$$B_k=\{i, 1\leq i \leq p,  \; s.t., a_i\in B(y_k, \rho_0)\}$$and $i_k$ an index of $B_k$ such that
\begin{equation}\label{eq1.0}
    \lambda_{i_k}=\min\{\lambda_i, i\in B_k\}.
\end{equation}
Denote $B_1, \ldots, B_\ell$ all the sets such that $\sharp B_k\geq 2, \forall k=1, \ldots, \ell$.

\n For any $\gamma>0$ and small, we define the following increasing cut-off function
\begin{eqnarray}
 \nonumber
 \chi: \mathbb{R}&\rightarrow&\mathbb{R}\\
 \nonumber   t&\mapsto&\left\{
             \begin{array}{ll}
               0, & \hbox{if } |t|\leq \gamma\\
               1/2, & \hbox{if } |t|=2\delta\\
               1& \hbox{if } |t|\geq 1.
             \end{array}
           \right.
\end{eqnarray}
For any $j\in B_k, k=1, \ldots, \ell$, we set $\bar{\chi}(\lambda_j)$ the following real number
$$\bar{\chi}(\lambda_j)= \sum_{i\in B_k, i\neq j}\chi\Big(\frac{\lambda_j}{\lambda_i}\Big),$$
and $Y_j(u)$ is the vector field defined in \eqref{2.B}.

\n We say that $j\in (R)$ for $j\in B_k, k=1, \ldots, \ell$, if there exists an index $i_j\in B_k$ such that $i_j\neq j$ and satisfying the following symmetric relation

\begin{equation}\label{R}
    2\gamma \lambda_j\leq\lambda_{i_j}\leq\frac{1}{2\gamma}\lambda_j.
\end{equation}
In the next, we denote

\begin{equation}\label{eq1.0'}
B_k^*= B_k, \text{ if } i_k \in (R) \text{ and } B_k^*= B_k\setminus\{i_k\}, \text{ if } i_k\not\in (R),
\end{equation}
for any $k=1, \ldots, \ell$. We now introduce the following pseudogradient
$$Z(u)= -\sum_{k=1}^\ell \sum_{j\in B_k^*} \bar{\chi}(\lambda_j) Z_j(u).$$

\begin{prop}\label{prop5}
There exists a bounded pseudogradient $\Xi_1$ on $V_3(p, \varepsilon)$ satisfying  $(P.S)$ condition and
$$<\partial J(u), \Xi_1(u)> \leq -c \sum_{k=1}^\ell\sum_{j\in B_k^*} \bar{\chi}(\lambda_j) \Big(\sum_{ j\neq i} \varepsilon_{ij}+ O\Big(\frac{1}{\lambda_j^{\beta_j}}\Big)\Big).$$
\end{prop}

\begin{pf}
Using expansion $(ii)$ of Proposition \ref{prop4}, we have
$$<\partial J(u), Z(u)>= 2J(u) \sum_{k=1}^\ell \sum_{j\in B_k^*}\bar{\chi}(\lambda_j)\Big[\tilde{c}_1 \sum_{j\neq i} \alpha_i \alpha_j \lambda_j \frac{\partial \varepsilon_{ij}}{\partial \lambda_j} + \sum_{s=2}^{[\beta]} O\Big(\frac{|a_j-y_k|^{\beta_j-s}}{\lambda_j^{s}}\Big)$$
\begin{equation}\label{eq1.1}
 +O\Big(\frac{1}{\lambda_i^{\beta_j}}\Big)+ o\Big(\sum_{j\neq i}\varepsilon_{ij}\Big)\Big].
\end{equation}
Let $j\in B_k^*, k=1, \ldots, \ell$. A direct computation shows that
\begin{equation}\label{eq1.2}
\lambda_j \frac{\partial \varepsilon_{ij}}{\partial\lambda_j}< -c \varepsilon_{ij},  \text{ if } i \not\in B_k,
\end{equation}
\begin{equation}\label{eq1.3}
\lambda_j \frac{\partial \varepsilon_{ij}}{\partial\lambda_j}\leq -c \varepsilon_{ij},  \text{ if } i \in B_k \text{ and } \lambda_i \text{ and } \lambda_j \text{ are of the same order,}
\end{equation}
\begin{equation}\label{eq1.4}
\bar{\chi}(\lambda_j) \lambda_j \frac{\partial \varepsilon_{ij}}{\partial\lambda_j} + \bar{\chi}(\lambda_i) \lambda_i \frac{\partial \varepsilon_{ij}}{\partial\lambda_i}\leq -c \varepsilon_{ij},  \text{ if } i \in B_k \text{ and } \lambda_i \text{ and } \lambda_j \text{ are not of the same order.}
\end{equation}

\begin{rem}
We notice that in \eqref{eq1.3}, $\lambda_i$ and $\lambda_j$ are of the same order means that $\dis \frac{1}{2} \lambda_i \leq \lambda_j\leq 2\lambda_i$. Therefore the parameter $c$ in \eqref{eq1.2}\ldots \eqref{eq1.4} is independent of $\gamma$.
\end{rem}

Let $M_0$ be a positive constant large enough so that if $\lambda_j |a_j-y_k|\geq M_0$, then
\begin{equation}\label{eq1.5}
\frac{1}{\lambda_j^{\beta_j}}=o\Big(\frac{|a_j-y_k|^{\beta_j-1}}{\lambda_j}\Big) \text{ and } \frac{|a_j-y_k|^{\beta_j-s}}{\lambda_j^s}= o\Big(\frac{|a_j-y_k|^{\beta_j-1}}{\lambda_j}\Big), s\geq 2.
\end{equation}
Let $\varphi_0$ be the following cut-off function
\begin{eqnarray}
 \nonumber
 \varphi_0: \mathbb{R}&\rightarrow&\mathbb{R}\\
 \nonumber   t&\mapsto&\left\{
             \begin{array}{ll}
               0, & \hbox{ if } |t|\leq M_0,\\
               1, & \hbox{ if } |t|\geq 2M_0.
             \end{array}
           \right.
\end{eqnarray}
Using \eqref{eq1.2}\ldots \eqref{eq1.5}, the estimate \eqref{eq1.1} is then reduced to
\begin{equation}\label{eq1.6}<\partial J(u), Z(u)>\leq -c \sum_{k=1}^\ell \sum_{j\in B_k^*}\bar{\chi}(\lambda_j)\Big[ \sum_{j\neq i} \varepsilon_{ij}  + O\Big(\frac{1}{\lambda_j^{\beta_j}}\Big)
+ o\Big(\varphi_0(\lambda_j |a_j-y_k|) \frac{|a_j-y_k|^{\beta_j-1}}{\lambda_j}\Big)\Big].
\end{equation}
Let $\Xi_1$ be the following vector field,
$$\Xi_1(u)= Z(u)+ \sum_{k=1}^\ell \sum_{j\in B_k^*} Y_j(u),$$
where $Y_j(u)$ is the vector field defined in \eqref{2.B}, where the cut-off function $\varphi_{M}$ equals to $\varphi_{0}$ and $m_0$ is a fixed positive constant small enough.

\n Using expansion $(i)'$ of Proposition \ref{prop4} and \eqref{eq1.5}, we have
\begin{equation}\label{eq1.7}
<\partial J(u), Y_j(u) >\leq -c \bar{\chi}(\lambda_j) \varphi_0(\lambda_j |a_j-y_k|) \Big[\frac{|a_j-y_k|^{\beta_j-1}}{\lambda_j}+ O\Big(\sum_{j\neq i} \varepsilon_{ij}\Big) \Big].
\end{equation}
Thus, for $m_0$ small, (depending only on the function $K$), we get from \eqref{eq1.6} and \eqref{eq1.7}
$$<\partial J(u), \Xi_1(u)> \leq -c \sum_{k=1}^\ell \sum_{j\in B_k^*} \bar{\chi}(\lambda_j) \Big[\sum_{j\neq i} \varepsilon_{ij} + O\Big(\frac{1}{\lambda_j^{\beta_j}}\Big)+ \frac{m_0}{2} \varphi_0(\lambda_j |a_j-y_k|) \frac{|a_j-y_k|^{\beta_j-1}}{\lambda_j} \Big]
$$$$ \leq -c \Big[\sum_{k=1}^\ell \sum_{j\in B_k^*} \bar{\chi}(\lambda_j) \Big(\sum_{j\neq i} \varepsilon_{ij} + O\Big(\frac{1}{\lambda_j^{\beta_j}}\Big) \Big].
$$
\end{pf}

\begin{prop}\label{prop6}
There exists a bounded pseudogradient $\Xi_2$ on $V_3(p, \varepsilon)$ satisfying  $(P.S)$ condition and the following estimate
$$<\partial J(u), \Xi_2(u)> \leq -c \Big[\sum_{k=1}^\ell\sum_{j\in B_k^*, j\in (R)} \bar{\chi}(\lambda_j) \Big(\sum_{ j\neq i} \varepsilon_{ij}+ \frac{1}{\lambda_j^{\beta_j}}\Big)+ \sum_{k=1}^\ell\sum_{j\in B_k^*, j\not\in (R)} \bar{\chi}\Big(\sum_{ j\neq i} \varepsilon_{ij}+ O\Big(\frac{1}{\lambda_j^{\beta_j}}\Big)\Big)\Big].$$
\end{prop}

\begin{pf}
Let $M=M(\gamma, m_0)$ be a positive constant larger than $2 M_0$. Here $m_0$ $M_0$ are fixed in the proof of Proposition \ref{prop5}.

\n For \begin{eqnarray}
 \nonumber
 \varphi: \mathbb{R}&\rightarrow&\mathbb{R}\\
 \nonumber   t&\mapsto&\left\{
             \begin{array}{ll}
               0, & \hbox{ if } |t|\leq M,\\
               1, & \hbox{ if } |t|\geq 2M.
             \end{array}
           \right.
\end{eqnarray}
We define the following pseudogradient,
$$Y(u)= \sum_{k=1}^\ell\sum_{j\in B_k^*, j\in (R)} Y_j(u),$$
where $Y_j(u)$ is the vector field defined in \eqref{2.B}, taking $\varphi_{M}=\varphi$. Since $M> 2M_0$, we get by \eqref{eq1.7} (we replace $\varphi_0$ by $\varphi$),
\begin{equation}\label{eq1.8}
<\partial J(u), Y(u)> \leq -c \sum_{k=1}^\ell\sum_{j\in B_k^*, j\in (R)} \bar{\chi}(\lambda_j)\varphi(\lambda_j|a_j-y_k|) \Big[\frac{|a_j-y_k|^{\beta_j-1}}{\lambda_j} + O\Big(\sum_{ j\neq i} \varepsilon_{ij}\Big)\Big].
\end{equation}
Let $$\Xi_2(u)= \Xi_1(u)+m_0 Y(u).$$
Using \eqref{eq1.8} and the estimate of Proposition \ref{prop5}, we have
$$<\partial J(u), \Xi_2(u)> \leq -c \Big[\sum_{k=1}^\ell\sum_{j\in B_k^*, j\in (R)} \bar{\chi}(\lambda_j) \Big(\sum_{ j\neq i} \varepsilon_{ij}+ O\Big(\frac{1}{\lambda_j^{\beta_j}}\Big)+ m_0 \varphi(\lambda_j|a_j-y_k|) \frac{|a_j-y_k|^{\beta_j-1}}{\lambda_j}\Big)$$
\begin{equation}\label{eq1.9}
 + \sum_{k=1}^\ell\sum_{j\in B_k^*, j\not\in (R)} \bar{\chi}\Big(\sum_{ j\neq i} \varepsilon_{ij}+ O\Big(\frac{1}{\lambda_j^{\beta_j}}\Big)\Big)\Big].
\end{equation}
For any $j_0\in B_k^*, k=1, \ldots, \ell$, such that $j_0\in (R)$, there exists $i_0 \in B_k$ such that $i_0\neq j_0$ and $$2\gamma\lambda_{j_0} \lambda_{i_0}< \frac{1}{2\gamma}\lambda_{j_0}.$$
Three cases may occur.

\n Case 1. If $\lambda_{i_0}|a_{i_0}-y_k|\geq 2 M$. We claim that
$$\frac{1}{\lambda_{j_0}^{\beta_{j_0}}}= o\Big(m_0 \frac{|a_{i_0}-y_k|^{\beta_{i_0}-1}}{\lambda_{i_0}}\Big), \text{ as  $M$ large.}$$
Indeed, using the fact that $\beta_{j_0}=\beta_{i_0}=\beta(y_k)$, we have
$$\frac{1}{m_0} \frac{1}{\lambda_{j_0}^{\beta_{j_0}}} \frac{\lambda_{i_0}}{|a_{i_0}-y_k|^{\beta_{i_0}-1}}\leq \frac{1}{m_04\gamma^{j_0}} \frac{1}{\lambda_{i_0}|a_{i_0}-y_k|^{\beta_{i_0}-1}}=o(1),$$as $M=M(\gamma, m_0)$ is large. Using the fact that $\bar{\chi}(\lambda_{i_0})\geq \frac{1}{2}$, $O\Big(\frac{1}{\lambda_{j_0}^{\beta_{j_0}}} \Big)$ which appears in the upper bound of \eqref{eq1.9} will be absorbed by $m_0 \bar{\chi}(\lambda_{i_0})\frac{|a_{i_0}-y_k|^{\beta_{i_0}-1}}{\lambda_{i_0}}$ and appears of the form $\frac{1}{\lambda_{j_0}^{\beta_{j_0}}}$.

\n Case 2. If $\lambda_{j_0}|a_{j_0}-y_k|\geq 2 M$. By \eqref{eq1.5}, we have
$$\frac{1}{\lambda_{j_0}^{\beta_{j_0}}}= o\Big(m_0 \frac{|a_{j_0}-y_k|^{\beta_{j_0}-1}}{\lambda_{j_0}}\Big), \text{ as  $M\geq 2M_0$.}$$
\n Case 3. If $\lambda_{i}|a_{i}-y_k|\leq 2 M$ for $i=i_0, j_0$. In this case, we claim that
$$\frac{1}{\lambda_{j_0}^{\beta_{j_0}}}= o\Big(\varepsilon_{i_0 j_0}\Big), \text{ as  $\varepsilon$ small.}$$
Indeed,

$$\frac{1}{\lambda_{j_0}^{\beta_{j_0}}} \varepsilon_{i_0 j_0}^{-1}= \bigg(\frac{\lambda_{i_0}}{\lambda_{j_0}^{\frac{2\beta_{j_0}}{n-2}+1}}+
\frac{1}{\lambda_{i_0}\lambda_{j_0}^{\frac{2\beta_{j_0}}{n-2}-1}}+ \frac{\lambda_{i_0}\lambda_{j_0}|a_{i_0}-a_{j_0}|^2 }{\lambda_{j_0}^{\frac{2\beta_{j_0}}{n-2}}}\bigg)^\frac{n-2}{2}$$
$$\leq \frac{1}{(2\gamma)^\frac{n-2}{2}} \frac{1}{\lambda_{j_0}^{\beta_{j_0}}}\Big(2+\lambda_{j_0}^2|a_{i_0}-a_{j_0}|^2 \Big)^\frac{n-2}{2}.$$
Observe that
\begin{eqnarray}
\nonumber
\lambda_{j_0}^2|a_{i_0}-a_{j_0}|^2 &\leq& \lambda_{j_0}^2|a_{i_0}-y_k|^2+ \lambda_{j_0}^2|a_{j_0}-y_k|^2\\ \nonumber
&\leq& \frac{1}{4\gamma^2}\lambda_{i_0}^2|a_{i_0}-y_k|^2+ \lambda_{j_0}^2|a_{j_0}-y_k|^2\\ \nonumber
&\leq& 2(\frac{1}{4\gamma^2} +1) M.
\end{eqnarray}
Thus, $\dis \frac{1}{\lambda_{j_0}^{\beta_{j_0}}} \varepsilon_{i_0 j_0}^{-1} \rightarrow 0$,  as  $\varepsilon \rightarrow 0$ and our claim follows.

\n In this case, $O\Big(\frac{1}{\lambda_{j_0}^{\beta_{j_0}}}\Big)$ will be absorbed  by $\varepsilon_{i_0 j_0}$, for $\varepsilon$ small enough and appeared of the form $\frac{1}{\lambda_{j_0}^{\beta_{j_0}}}$. This conclude the proof of Proposition \ref{prop6}.
\end{pf}

\begin{cor}\label{cor3}
For $\gamma>0$ small enough, we have
$$<\partial J(u), \Xi_2(u)> \leq -c \bigg[\sum_{k=1, i_k\in (R)}^\ell\sum_{j\in B_k} \bar{\chi}(\lambda_j) \Big(\sum_{ j\neq i} \varepsilon_{ij}+ \frac{1}{\lambda_j^{\beta_j}}\Big)$$$$+ \sum_{k=1, i_k\not\in (R)}^\ell \bigg(\sum_{j\in B_k^*, j\in (R)} \bar{\chi}(\lambda_j)\Big(\sum_{ j\neq i} \varepsilon_{ij}+ \frac{1}{\lambda_j^{\beta_j}}\Big)+ \sum_{j\in B_k^*, j\not \in (R)} \bar{\chi}(\lambda_j)\Big(\sum_{ j\neq i} \varepsilon_{ij}+ O\Big(\frac{1}{\lambda_j^{\beta_j}}\Big)\Big)\bigg)\bigg].$$
\end{cor}

\begin{pf}
For any $k_0=1, \ldots, \ell$, such that $i_{k_0}\in (R)$, we have $B_{k_0}^*= B_{k_0}$. For any $j\in B_{k_0}$ such that $j\not\in (R)$, we have
$$\lambda_{i_{k_0}}< 2\gamma \lambda_j.$$
Thus, $\frac{1}{\lambda_j^{\beta_j}}=o\Big(\frac{1}{\lambda_{i_{k_0}}^{\beta_{i_{k_0}}}}\Big)$, as $\gamma$ small. Therefore, $O\Big(\frac{1}{\lambda_j^{\beta_j}}\Big)$ which appears in the upper  bound of the estimate of Proposition \ref{prop6}, will be absorbed by $\bar{\chi}(\lambda_{i_k})\frac{1}{\lambda_{i_{k}}^{\beta_{i_{k}}}}$, ($\bar{\chi}(\lambda_{i_k})\geq \frac{1}{2}$) and appears of the form $\frac{1}{\lambda_j^{\beta_j}}$.
\end{pf}

\n Now let $$D=\{i_k, i_k=1 \ldots, \ell, s.t., i_k\not\in (R)\}\cup \Big(\dis \cup_{k=1}^\ell B_k\Big)^c.$$
For $\hat{u}= \sum_{i\in D}\alpha_i \delta_{(a_i, \lambda_i)}$, we have $y_i\neq y_j, \forall i\neq j\in D$. Therefore, $\hat{u}\in W_1(\sharp D, \varepsilon)$ or $\hat{u}\in W_2(\sharp D, \varepsilon)$. We apply the corresponding pseudogradient $\widetilde{W}_i(\hat{u}), i=1$ or  $2$. We have

\begin{equation}\label{eq1.10}
<\partial J(u), \widetilde{W}_i(u)> \leq -c \bigg(\sum_{j\in D} \Big(\frac{1}{\lambda_j^{\beta_j}}+\frac{|\nabla K(a_j)|}{\lambda_j}\Big) + \sum_{ i\neq j\in D} \varepsilon_{ij}\bigg)+ \sum_{j\in D, i\not \in D}O( \varepsilon_{ij}).
\end{equation}
Let $W_3(u)= \Xi_2(u)+ m_0 \widetilde{W}_i(\hat{u})$. Using \eqref{eq1.10} and the result of Corollary \ref{cor3}, we have
\begin{equation}\label{eq1.11}
<\partial J(u), {W}_3(u)> \leq -c \bigg(\sum_{i=1}^p \frac{1}{\lambda_i^{\beta_i}}+\sum_{i\in D}\frac{|\nabla K(a_i)|}{\lambda_i} + \sum_{ i\neq j} \varepsilon_{ij}\bigg),
\end{equation}
since for $j\in B_k^*$, s.t. $i_k\not \in (R)$, and $j\not \in (R)$, we have $$\frac{1}{\lambda_j^{\beta_j}}= o \Big(\frac{1}{\lambda_{i_k}^{\beta_{i_k}}}\Big), \text{ as $\gamma$ small enough.}$$
In order to make appear $\dis -\sum_{j=1}^p\frac{|\nabla K(a_j)|}{\lambda_j}$ in the upper bound of \eqref{eq1.11}, we consider
$$\widetilde{W}_3(u)= W_3(u)+ m_0\sum_{j=1, j\not \in D}^p Y_j(u).$$
It satisfies
$$<\partial J(u), \widetilde{W}_3({u})>\leq -c \bigg( \sum_{i=1}^p\Big( \frac{1}{\lambda_{i}^{\beta_{i}}}+\frac{|\nabla K(a_i)|}{\lambda_i}\Big)+ \sum_{i\neq j}\varepsilon_{ij}\bigg).$$\\

\n $\bullet$ \underline{Pseudogradient in $V_4(p, \varepsilon), p\geq 1$:}
\vskip0.15cm

\n Let $u=\sum_{i=1}^{p}\alpha_i \delta_{(a_i, \lambda_i)}\in V_4(p, \varepsilon)$. Assume that
$$\lambda_1\leq ...\leq \lambda_p.$$
Let $j_1$ be the first index such that $a_{j_1}\not\in B(y,\rho_0), \ \forall y \in \Gamma$ and let $j_0$ be the first index such that $\lambda_{j_{0}}\geq \frac{1}{2}\lambda_{j_{1}}$. For a small positive constant $m$ we define\\
$$X(u)=\dis \frac{m}{\lambda_{j_{1}}}\dis\frac{\partial \delta_{(a_{j_{1}}, \lambda_{j_{1}})}}{\partial a_{j_{1}}}\dis\frac{\nabla K(a_{j_{1}})}{|\nabla K(a_{j_{1}})|}-\dis\sum_{i\geq j_0}2^{i}Z_i(u).$$\\
Since $|\nabla K(a_{j_{1}})|\geq c >0$, we get
$$<\partial J(u), X(u)>\leq - c \bigg(\sum_{i\geq j_0,i \neq j} \varepsilon_{ij} + \frac{1}{\lambda_{j_{1}}}\bigg)$$
$$\leq - c \sum_{i\geq j_{0}}\bigg( \frac{1}{\lambda_{i}^{\beta_{i}}}+ \frac{|\nabla K(a_i)|}{\lambda_{i}}+ \sum_{j \neq i} \varepsilon_{ij}\bigg).$$

\n Let $\widehat {u}=\sum_{i<j_{0}}\alpha_i \delta_{(a_i, \lambda_i)}$. $\widehat {u}$ has to satisfy the conditions of one of the above regions $V_i(j_{0}-1, \varepsilon), \ i=1,2,3$. Let $\widetilde {W}_i(\widehat {u})$ the corresponding vector field. For $m'>0$ and small, we set
$$\widetilde {W}_4(u)=X(u)+m'\widetilde {W}_i(\widehat {u}).$$
$\widetilde {W}_4$ satisfies $(P.S)$ condition, moreover we have
$$<\partial J(u), \widetilde{W}_4({u})>\leq -c \bigg( \sum_{i=1}^p\Big( \frac{1}{\lambda_{i}^{\beta_{i}}}+\frac{|\nabla K(a_i)|}{\lambda_i}\Big)+ \sum_{i\neq j}\varepsilon_{ij}\bigg).$$\\
Finally, let $W$ be a convex combination of $\widetilde{W}_i, \ i=1,...,4$. Define\\
$$\widehat{W}_1(u)=W(u)-\Big\langle W(u),u\Big\rangle_{H^1}u, \ u\in V(p, \varepsilon)$$
and
$$\widehat{W}_2(u)=-\partial J(u), \ u\in V(p, \frac{\varepsilon}{2})^c.$$
The global pseudogradient $\widetilde{W}$ in the variational space $\Sigma$ is defined by a convex combination of $\widehat{W}_1$ and $\widehat{W}_2$. By construction $\widetilde{W}$ is bounded and satisfies $(a)$ and $(c)$ of Theorem \ref{th3.1}. Concerning claim $(b)$ it follows as in (\cite{bcch1}, Appendix 2) from $(a)$ and the estimate of $\|\bar{v}\|^2$ of Proposition \ref{prop3}. The  proof of Theorem \ref{th3.1} is thereby completed.
\end{pfn}

\vspace{0.1cm}
\n Let $\varepsilon_0$ be a fixed positive constant small enough and let\\
 \begin{center}
 $V_{\varepsilon_0}(\Sigma^+)=\big\{u \in \Sigma,  \, J(u)^{\frac{n}{2}}\|u^{-}\|<\varepsilon_0 \big\}$.
\end{center}
\begin{lem}$V_{\varepsilon_0}(\Sigma^+)$ is an invariant set under the action of the flow of $\widetilde {W}$.
\end{lem}
\begin{pf} Since $\widetilde {W}=-\partial J$ in $A$, where
$$A=\big\{u \in \Sigma,  \, \frac{\varepsilon_0}{2}<J(u)^{\frac{n}{2}}\|u^{-}\|<\varepsilon_0 \big\},$$
the proof proceeds exactly as the one of (\cite{bcch1}, Lemma 4.1).
\end{pf}

\n We now introduce the following main result.
\begin{thm}\label{th3.3}
Assume that \eqref{1.1} has no solution. Under assumptions $(f)_{\beta}$, $\mathbf{(H_1)}$and $\mathbf{(H_2)}$, the critical points at infinity of the variational functional associated to problem \eqref{1.1} are:
$$(y_1)_\infty:= \frac{1}{K(y_1)^\frac{n-2}{2}}\delta_{(y_1, \infty)}, \; y_1\in \Gamma^-,$$and
$$(y_1, \ldots, y_p)_\infty:= \sum_{i=1}^p\frac{1}{K(y_i)^\frac{n-2}{2}}\delta_{(y_i, \infty)}, \; y_i\neq y_j, \ \forall i\neq j, \hbox{ and }\{y_1, \ldots, y_p\}\in \widetilde{\Lambda}^-.$$
The index of each critical point at infinity $(y_1, \ldots, y_p)_\infty$ equals to $p-1+\sum_{j=1}^p n-i(y_j).$
\end{thm}
\vspace{0.1cm}
\textbf{Proof of Theorem \ref{th3.3}}  Assume that \eqref{1.1} has no solution. Let $u_0$ be an initial condition in $V_{\varepsilon_0}(\Sigma^+)$. Since in $\biggr[\dis\bigcup_{p}V(p, \frac{\varepsilon}{2})\biggr]^c$ we have $<\partial J(u), \widetilde{W}(u)>\leq-c$, then there exist $s_0=s(u_0)>0$ and $p=p(u_0)\geq 1$ such that the flow line $\eta(s,u_0)$ generated by $\widetilde{W}$ lies in $V(p, \frac{\varepsilon}{2})$, for any $s\geq s_0$. As seen, the pseudogradient $\widetilde{W}$ built in Theorem \ref{th3.1} satisfies $(P.S)$ condition on its flow lines as long as these flow lines do not enter a small neighborhood $\mathcal{N}(y_1, \ldots, y_p)$ of $\dis \sum_{i=1}^p \frac{1}{K(y_i)^{\frac{n-2}{2}}}\delta_{(y_i, \infty)}$ such that $y_1 \in \Gamma^-$ if $p=1$ and $y_i\neq y_j, \ \forall 1\leq i\neq j \leq p \hbox{ with }\{y_1, \ldots, y_p\}\in \widetilde{\Lambda}^-$ if $p\geq 2$. However, if a flow line of $\widetilde{W}$ enter in $\mathcal{N}(y_1, \ldots, y_p)$, all the concentrations $\lambda_i, i=1, \ldots, p$ tend to $\infty$ and the flow remains in $\mathcal{N}(y_1, \ldots, y_p)$ for all time $s\geq s_1$. This conclude the characterization of critical points at infinity.

\n Arguing as in (\cite{bcch1}, proof of Lemma 4.2), the functional $J$ can be expanded in $\mathcal{N}(y_1, \ldots, y_p)$ as follows

$$J\bigg(\sum_{i=1}^p \alpha_i \delta_{(a_i, \lambda_i)}+ \bar{v}\bigg) = S_n \bigg(\sum_{i=1}^p \frac{1}{K(y_i)^{\frac{n-2}{2}}}\bigg)^{\frac{2}{n}}\bigg(1-|H|^2+ \sum_{i=1}^p \Big(|a_i^-|^2 - |a_i^+|^2\Big)\bigg),$$
where $H\in \mathbb{R}^{p-1}$ and $(a_i^-, a_i^+)$ denote the coordinate of $a_i, i=1, \ldots, p$ along the unstable manifold and the stable manifold of $K$ at $y_i$. We then derive that the index of $J$ at $\dis \sum_{i=1}^p \frac{1}{K(y_i)^{\frac{n-2}{2}}}\delta_{(y_i, \lambda_i)}$ is equal to $p-1+ \sum_{i=1}^p (n-\widetilde{i}(y_i))$. This finishes the proof of Theorem \ref{th3.3}. $\square$

\section{proof of results}

\n We now prove the Theorems of the first section. We need to introduce the following lemma
\begin{lem}\label{lem3.13} Let $w$ be a solution of \eqref{1.1}. Under assumptions $(f)_{\beta}$ and $\mathbf{(A_1)}$, $V(p,\varepsilon,w)$ contains no critical point at infinity for any $p\geq 1$.
\end{lem}
\begin{pf}Let $u=\dis\sum_{i=1}^{p}\alpha_{i}
 \delta_{(a_{i},\lambda_{i})}+\alpha_{0}
(\omega+h)+\overline{v} \in V(p,\varepsilon,w)$. Following (\cite{bch}, Proof of Theorem 1.3), there exists a change of
variables
$$(a_{i},\lambda_{i},h)\longmapsto(\tilde a_{i},\tilde\lambda_{i},\tilde h),$$
 such that $$J\Big(\dis\sum_{i=1}^{p}\alpha_{i}
 \delta_{(a_{i},\lambda_{i})}+\alpha_{0}
(\omega+h)+\overline{v}\Big)=J\Big(\dis\sum_{i=1}^{p}\alpha_{i}\delta_{(\tilde a_{i},\tilde\lambda_{i})}
+\alpha_{0}(w+\tilde h)\Big).$$
By extending the computation of (\cite{BOA1}, Proposition 3.1) to the case of $(f)_{\beta}$ condition, we have
$$J\big(\dis\sum_{i=1}^{p}\alpha_{i}
 \delta_{(a_{i},\lambda_{i})}+\alpha_{0}
(\omega+h)\big)$$
\begin{eqnarray*}
&=&\dis\frac{S_{n}\dis\sum_{i=1}^{p}\alpha_{i}^{2}+\alpha_{0}^{2}\parallel
w\parallel^{2}}{(S_{n}\dis\sum_{i=1}^{p}\alpha_{i}^{\frac{2n}{n-2}}K(a_{i})+\alpha_{0}^{\frac{2n}{n-2}}\parallel
w\parallel^{2})^{\frac{n-2}{n}}}
\biggl\{1-c\alpha_{0}\dis\sum_{i=1}^{p}
\alpha_{i}\dis\frac{w(a_{i})}{\lambda_{i}
^{\frac{n-2}{2}}}-c'\sum_{i=1}^p\alpha_i^{\frac{2n}{n-2}}\frac{\sum_{k=1}^n b_k(y_i)}{\lambda_i^{\beta(y_i)}}\\
&&-c''\dis\sum_{i\neq
j}\alpha_{i}\alpha_{j}\varepsilon_{ij}
\biggr\}+Q(h,h),
\end{eqnarray*}
provided $a_i \in B(y_i,\rho)$, $y_i\in \Gamma, \ \forall i=1,...,p$. Here $Q(h,h)$ is a negative quadratic form, (see \cite{BOA1}, Lemma 3.2).\\
Therefore, if $\beta >\frac{n-2}{2}$, the contribution of $w$ in the above expansion becomes dominant with respect to the term $\frac{1}{\lambda_i^{\beta_i}}$ and pushes  down any flow line in $V(p,\varepsilon,w)$.
\end{pf}

\n\textbf{Proof of Theorem \ref{th1.2}} Using (\cite{yy1}, section 6) and (\cite{L3}, Theorem 09), the proof of Theorem \ref{th1.2} follows from Theorem \ref{th3.3} and Lemma \ref{lem3.13}.$\square$ \\

\n\textbf{Proof of Theorem \ref{th1.1}} It follows from Theorem \ref{th1.2}.$\square$ \\

\n\textbf{Proof of Theorem \ref{th1.3}} We use $\widetilde{W}$ to deform $V_{\varepsilon_0}(\Sigma^+)$. It follows from the results of Theorem \ref{th3.3} and Lemma \ref{lem3.13} that under assumptions of Theorem \ref{th1.3}, the only critical points at infinity of $J$ are
$$\frac{1}{K(y)^\frac{n-2}{2}}\delta_{(y, \infty)}, \; y\in \Gamma^-,$$and
$$\sum_{i=1}^p\frac{1}{K(y_i)^\frac{n-2}{2}}\delta_{(y_i, \infty)}, \; y_i\neq y_j, \ \forall i\neq j, \hbox{ with }\{y_1, \ldots, y_p\}\in \widetilde{\Lambda}^-.$$
Using the deformation Lemma of \cite{BR1}, we have

\begin{equation}\label{4.1}
V_{\varepsilon_0}(\Sigma^+)\simeq \bigcup_{y\in \Gamma^{-}} W_u^{\infty}
(y)_\infty \bigcup \bigcup_{A\in \widetilde{\Lambda}^{-}}W_u^\infty(A)_\infty\bigcup \bigcup_{w\in V_{\varepsilon_0}(\Sigma^+), \partial J(w)=0 }W_u(w).
\end{equation}
 Here $\simeq $ denotes retract by deformation , $W_u^\infty$ denotes the unstable manifold of a critical point at infinity and $W_u$ denotes the unstable manifold of a critical point of $J$.\\
By Sard-Smale Theorem \cite{ss}, we know for generic $K$, all critical points of $J$ can be considered as non degenerate critical points. By applying the Euler-Poincar\'{e} characteristic on the both side of \eqref{4.1}, we obtain
\begin{equation}\label{bb}
1= \sum_{y\in \Gamma^{-}} (-1)^{n-i(y)}+ \sum_{A=\{y_1, \ldots, y_p\}\in \widetilde{\Lambda}^{-}} (-1)^{p-1 + \sum_{j=1}^{p} n-i(y_j)}+ \sum_{w\in V_{\varepsilon_0}(\Sigma^+), \partial J(w)=0} (-1)^{i(w)}.
\end{equation}
From the above equality we deduce first that the functional $J$ has at least a critical point $u_0$ in $V_{\varepsilon_0}(\Sigma^+)$. Indeed,\\
if not the degree $d$ of Theorem \ref{th1.2} equals to zero which is a contradiction.\\
Arguing as the proof of the Theorem of (\cite{bcch1}, Pages 659-660) we can prove that $u_0^-=0$ and therefore $u_0$ is a solution of \eqref{1.1}.\\
Concerning the number of solutions, it follows from \eqref{bb} that
$$|d|=\big|\sum_{w\in V_{\varepsilon_0}(\Sigma^+), \partial J(w)=0} (-1)^{i(w)}\big|.$$
This finishes the proof of Theorem \ref{th1.3}.$\square$ \\

\n\textbf{Proof of Theorem \ref{th1.4}} It follows from Theorem \ref{th3.3} and the argument of the proof of Theorem \ref{th1.3}.$\square$

\end{document}